# OPTIMAL CONSUMPTION IN DISCRETE-TIME FINANCIAL MODELS WITH INDUSTRIAL INVESTMENT OPPORTUNITIES AND NONLINEAR RETURNS

By Bruno Bouchard and Huyên Pham

*Université Paris 6 and Université Paris 7*

We consider a general discrete-time financial market with proportional transaction costs as in [Kabanov, Stricker and Rásonyi *Finance and Stochastics* **7** (2003) 403–411] and [Schachermayer *Math. Finance* **14** (2004) 19–48]. In addition to the usual investment in financial assets, we assume that the agents can invest part of their wealth in industrial projects that yield a nonlinear random return. We study the problem of maximizing the utility of consumption on a finite time period. The main difficulty comes from the nonlinearity of the nonfinancial assets' return. Our main result is to show that existence holds in the utility maximization problem. As an intermediary step, we prove the closedness of the set $A_T$ of attainable claims under a *robust no-arbitrage* property similar to the one introduced in [Schachermayer *Math. Finance* **14** (2004) 19–48] and further discussed in [Kabanov, Stricker and Rásonyi *Finance and Stochastics* **7** (2003) 403–411]. This allows us to provide a dual formulation for $A_T$.

**1. Introduction.** We consider a general discrete-time market with proportional transaction costs as in [7, 8, 13]. Following the above papers, we model the wealth process by a vector-valued process $(V_t)$, each component $i$ corresponding to the number of units of asset $i$ which are held in the portfolio. The usual self-financing condition is described by the constraints $V_t - V_{t-1} \in -K_t$, where $-K_t$ is the random convex set of affordable exchanges at time $t$, given the value of the underlying assets and the level of transaction costs.

In the case of efficient frictions where the transaction costs are positive (which is formulated by the assumption that $K_t$ is proper), a general version

Received September 2004; revised January 2005.

of the Fundamental Theorem of Asset Pricing was obtained by Kabanov, Stricker and Rásonyi [7]. In the case where some of the costs may be zero, a notion of "robust no-arbitrage" was introduced by Schachermayer [13] and further studied in [8]. This assumption can be interpreted as follows: there is no arbitrage even if we reduce the size of the proportional transaction costs (which are not already equal to zero). In the above papers, it is shown that this assumption is equivalent to the existence of a strictly consistent price system (see [13] for a precise definition). It also implies the closedness of the set of attainable claims and allows to provide a suitable dual formulation for this set.

In addition to the above setting, we assume in this paper that the financial agent can invest part of its wealth in nonfinancial assets, for example, industrial projects, which are also subject to proportional costs (see [4, 5]), but, in opposition to usual financial assets, yield nonlinear returns. Our principal aim is to study the problem of maximizing the utility of consumption over a finite time period. The analysis of such a model differs from the usual setting in many aspects:

1. It follows from the nonlinearity of the nonfinancial assets' return that the set $A_T(0)$ of attainable claims with zero initial endowment is not a cone. More generally, the set of attainable claims with initial endowment $x$, $A_T(x)$, is not linear with respect to $x$, that is, $x + A_T(0) \neq A_T(x)$.
2. All transactions $V_t - V_{t-1} \in -K_t$ are not allowed since it is natural to impose a nonnegativity constraint on the level of investment in the nonfinancial assets. In fact, the effective set of possible transactions at time $t$ is a subset of $-K_t$ which depends on the initial endowment and all the transactions up to $t$.
3. The notion of no-arbitrage is not as clear as in pure financial market. Indeed, if we have an initial investment $y$ (in units) in some project which yields a nonnegative return in terms of cash, and if we do nothing, at the time horizon $T$ we end up with a nonnegative amount of cash $g$ and we still have the investment $y$ (in units). Since $(g, y) \geq (0, y)$ there is an arbitrage, in the usual sense, if $g \neq 0$. However, from an economic point of view this situation should be possible as the risk supported by investing in a project also lies in the liquidation value of the investment which does not appear in the above formulation.

In order to avoid trivial situations, we have to impose some no-arbitrage condition. In view of part 3 above, we define it only on $A_T(0)$, that is, we assume that $A_T(0) \cap L^0(\mathbb{R}_+^{d+N}) = \{0\}$; see the notation below. As the initial endowment in nonfinancial asset is 0, this avoids the problem pointed out in part 3. In order to obtain the usual closedness property of $A_T(0)$, we impose a "robust no-arbitrage condition." Because of the nonlinearity of the nonfinancial assets' returns, we cannot work directly with the "robust no-arbitrage condition" of Schachermayer [13]. We therefore extend this definition. Our version can be interpreted as follows: there is no arbitrage even if we slightly reduce the size of the proportional transaction costs between financial assets and slightly increase the return of the nonfinancial ones. It also allows us to provide a dual formulation for this set.

In the multivariate setting, the usual duality approach for the utility maximization problem is much more complex than in the case of no transaction costs. The reason is that, even when the utility function $U$ is smooth (which is not assumed here), its Fenchel transform $\tilde{U}$ may not be smooth. To manage this difficulty, we can proceed as in Deelstra, Pham and Touzi [2] and Bouchard, Touzi and Zeghal [1] who reduce to the smooth case by approximating $\tilde{U}$ by smooth convex functions. But this leads to long and technical proofs. In the paper by Kramkov and Schachermayer [10], a more direct argument is proposed. It consists in first deriving the duality theorem in an abstract way. This allows to show that maximizing sequences for the primal problem satisfy a uniform integrability condition. However, it turns out that the one-dimensional argument of Kramkov and Schachermayer [10] does not work directly in our multivariate setting. We overcome this difficulty by introducing some auxiliary primal problem. At this point, we shall mention that another possible approach would be to use a dynamic programming as done in [11] for the (frictionless) incomplete discrete-time market. We leave this alternative for future research.

The rest of the paper is organized as follows. The model is described in Section 2. We discuss our "robust no-arbitrage" condition in Section 3. The utility maximization problem is defined in Section 4 where we state our existence result. In Section 5, we show the closedness of the set attainable terminal wealth and we provide a dual formulation for this set in Section 6. The last section contains the proof of the existence result.

In this paper we shall repeatedly use the following notation. For $x \in \mathbb{R}^{d+N}$, we shall often write $x$ as $(x^F, x^I)$ where $x^F \in \mathbb{R}^d$ and $x^I \in \mathbb{R}^N$. The exponent $F$ (resp. $I$) stands for "financial" (resp. "industrial"). Given $E \subset \mathbb{R}^{d+N}$, we write $\underline{E} = \{(x^F, 0_N) : x = (x^F, x^I) \in E\}$, where $0_N$ denotes the zero of $\mathbb{R}^N$. We denote by $\|\cdot\|$ the Euclidean norm and by "·" the inner product of $\mathbb{R}^p$, where $p \in \mathbb{N}$ is given by the context. $\mathbb{R}_+^p$ will denote the set of elements of $\mathbb{R}^p$ with nonnegative components. Given a probability space $(\Omega, \mathcal{F}, \mathbb{P})$ endowed with a filtration $\mathbb{F} = (\mathcal{F}_t)_{t \in \mathbb{T}}$, $\mathbb{T} = \{0, \ldots, T\}$ for some $T \in \mathbb{N} \setminus \{0\}$, and a random set $E$, we denote by $L^0(\Omega \times \mathbb{T}, E)$ the set of processes $Y = (Y_t(\omega))_{t \in \mathbb{T}}$ valued in $E$, by $L^0(E; \mathcal{F}_t)$ the set of $\mathcal{F}_t$-measurable random variables which take values in $E$ $\mathbb{P}$-a.s. For $\mathbb{F}$-adapted processes with values in $E$, we write $L^0(E; \mathbb{F})$. For $\tilde{\mathbb{P}} \sim \mathbb{P}$, we similarly denote by $L^1(\Omega \times \mathbb{T}, \tilde{\mathbb{P}}, E)$ [resp. $L^1(E; \tilde{\mathbb{P}}, \mathcal{F}_t)$] the set of elements of $L^0(\Omega \times \mathbb{T}, E)$ [resp. $L^0(E; \mathcal{F}_t)$] which are $\tilde{\mathbb{P}}$-integrable. For bounded random processes (resp. $\mathcal{F}_t$ measurable random variables), we use the notation $L^\infty(\Omega \times \mathbb{T}, E)$ [resp.

$L^\infty(E;\mathcal{F}_t)]$. When $\tilde{\mathbb{P}} = \mathbb{P}$, we omit the argument $\mathbb{P}$, and similarly when $t = T$, we may omit the argument $\mathcal{F}_t$. We do the same thing for $E$ when it is clearly given by the context. For a subset $E \subset \mathbb{R}^p$, we denote by $E^*$ its positive polar in the sense of convex analysis, that is, $E^* := \{y \in \mathbb{R}^p : x \cdot y \geq 0 \text{ for all } x \in E\}$. Given an event set $B$, we denote $E\mathbb{1}_B = \{\mathbb{1}_B x : x \in E\}$ where $\mathbb{1}_B = 1$ on $B$ and 0 otherwise. These last notations are naturally extended to random sets.

## 2. A financial model with industrial investment opportunities.

2.1. *Financial and industrial investment strategies.* Set $\mathbb{T} = \{0, \ldots, T\}$ for some $T \in \mathbb{N} \setminus \{0\}$ and let $(\Omega, \mathcal{F}, \mathbb{P})$ be a probability space endowed with a filtration $\mathbb{F} = (\mathcal{F}_t)_{t \in \mathbb{T}}$. We assume that $\mathcal{F}_T = \mathcal{F}$ and that $\mathcal{F}_0$ is trivial. Given two integers $d \geq 1$ and $N \geq 1$, we denote by $\mathcal{K}$ the set of $\mathcal{C}$-valued processes $K$ such that $\mathbb{R}^{d+N}_+ \setminus \{0\} \subset \text{int}(K_t)$ $\mathbb{P}$-a.s. for all $t \in \mathbb{T}$. Here, we follow Kabanov, Stricker and Rásonyi [8] and say that a sequence of set-valued mappings $(K_t)_{t \in \mathbb{T}}$ is a $\mathcal{C}$-valued process if there is a countable sequence of $\mathbb{R}^{d+N}$-valued processes $X^n = (X^n_t)_{t \in \mathbb{T}}$ such that for every $t \in \mathbb{T}$, $\mathbb{P}$-a.s. only a finite but nonzero number of $X^n_t$ is different from zero and $K_t = \text{cone}\{X^n_t, n \in \mathbb{N}\}$. This means that $K_t$ is the polyhedral cone generated by the $\mathbb{P}$-a.s. finite set $\{X^n_t, n \in \mathbb{N} \text{ and } X^n_t \neq 0\}$.

Given $K \in \mathcal{K}$, we denote by $\mathcal{A}(K)$ the set of processes $\xi \in L^0(\mathbb{R}^{d+N}; \mathbb{F})$ such that

$$\xi_t \in -K_t \quad \text{and} \quad I(\xi)^i_t := \sum_{s=0}^{t} \xi^{d+i}_s \geq 0, \qquad 1 \leq i \leq N, \mathbb{P}\text{-a.s. for all } t \in \mathbb{T}.$$

The interpretation is the following. For $1 \leq i \leq d + N$, the quantity $(\xi_t)^i$ corresponds to the number of units of asset $i$ which are bought at time $t$. The convex cone $-K_t$ is the set of variations in the global portfolio which are affordable, after possibly throwing out some units of the assets, at time $t$ given the price of the assets. The process $I(\xi)$ corresponds to the global investment in the different industrial projects. The condition $I(\xi)_t \in \mathbb{R}^N_+$, $\mathbb{P}$-a.s. means that it is not possible to have a negative level of investment in an industrial project.

Due to the constraint on the level of investment, we also need to consider the case where the strategy starts with an initial holding $x = (x^F, x^I) \in \mathbb{R}^d \times \mathbb{R}^N_+$. We then extend the previous notation and define $\mathcal{A}(x; K)$ as the set of processes $\xi \in L^0(\mathbb{R}^{d+N}; \mathbb{F})$ such that

(2.1) $\quad \xi_t \in -K_t \quad \text{and} \quad I(\xi)_t + x^I \in \mathbb{R}^N_+, \qquad \mathbb{P}\text{-a.s. for all } t \in \mathbb{T}.$

Observe that $\mathcal{A}(K) = \mathcal{A}((x^F, 0_N); K)$.

The return associated to the industrial investment is modeled by a process $R \in \mathcal{R}$, the set of adapted processes with values in the set of mapping from $\mathbb{R}_+^N$ into $\underline{\mathbb{R}}^{d+N}$; that is, $R_t(x)$ is $\mathcal{F}_t$-measurable for all $x \in \mathbb{R}_+^N$. A level of investment $I(\xi)_t$ in the industrial project at time $t$ leads to a reward (in units) $R_{t+1}(I(\xi)_t)$ at time $t+1$. Here, the fact that $R_{t+1}$ takes values in $\underline{\mathbb{R}}^{d+N}$ means that the reward consists in units of the financial assets. If the $N$ last assets are interpreted as industrial tools used for an industrial project, it is natural to assume that the reward consists in stocks or currencies, that is, pure financial assets, while the (relative) value of these tools may evolve in time.

The set of claims that can be reached with an initial holding $x = (x^F, x^I) \in \mathbb{R}^d \times \mathbb{R}_+^N$ is then given by

$$A_T(x; K, R) := \left\{ x + \sum_{t=0}^{T} \xi_t + \sum_{t=0}^{T-1} R_{t+1}(x^I + I(\xi)_t), \xi \in \mathcal{A}(x; K) \right\}.$$

For $x = 0$, we shall simply write $A_T(K, R)$ for $A_T(x; K, R)$.

REMARK 2.1. Observe from (2.1) that for general $x = (x^F, x^I) \in \mathbb{R}^d \times \mathbb{R}_+^N$, we do not have equality between $\mathcal{A}(x; K)$ and $\mathcal{A}(K)$, except if $x^I = 0$. Similarly, $A_T(x; K, R)$ differs from $x + A_T(K, R)$ in general, while $A_T(x; K, R) = x^F + A_T((0_d, x^I); K, R)$. Also, observe that $A_T(K, R)$ is in general not a cone since $R_t$ is not assumed to be linear.

In this paper, we shall assume that $(K, R) \in \mathcal{K} \times \mathcal{R}$ satisfies the above assumptions $\mathbb{P}$-a.s. for each $t \in \mathbb{T}$:

(R1) $R_t(0) = 0$ and $R_t$ is continuous.
(R2) For $\lambda \in [0, 1]$ and $(\alpha, \beta) \in (L^0(\mathbb{R}_+^N))^2$, we have

$$\lambda R_t(\alpha) + (1 - \lambda) R_t(\beta) - R_t(\lambda \alpha + (1 - \lambda) \beta) \in -\underline{K}_t.$$

(R3) There is some $a_t \in L^0(\mathbb{R}^{d+N})$ and $L \in \mathcal{R}$ such that $\lambda L_t(\alpha) = L_t(\lambda \alpha)$ $\mathbb{P}$-a.s. and $R_t(\alpha) + a_t + L_t(\alpha) \in L^0(\mathbb{R}_+^{d+N})$ for all $(\lambda, \alpha) \in L^0(\mathbb{R}_+ \times \mathbb{R}_+^N)$.

The condition $R_t(0) = 0$ is natural since no investment in the industrial project should yield no return. The condition (R2) is a concavity assumption. It means that, up to an immediate transaction in terms of financial assets, the return induced by a convex combination of industrial investments is better than the convex combination of the returns induced by each of them. It implies that $A_T(x; K, R)$ is convex (see Lemma 2.1 below). The last assumption is more technical. It imposes an affine lower bound on the mapping $x \mapsto R_t(x)(\omega)$ for almost every $\omega \in \Omega$. In the one-dimensional case, this means that $R'(\infty) > -\infty$ $\mathbb{P}$-a.s. It is used only in the proof of Lemma 5.3 below and can be replaced by a weaker one as explained in Remark 5.1.

Observe that we do not impose nonnegative returns; that is, an investment in nonfinancial assets may lead to a negative reward in terms of financial assets.

We conclude this section with a popular example in the economic literature.

EXAMPLE 2.1. Set $d=1$, $N=2$ and consider the Cobb–Douglas function $r(x,y) = x^\gamma y^\beta$ with $\beta$, $\gamma > 0$ and $\beta + \gamma < 1$. For $p \in L^0(\mathbb{R}_+; \mathbb{F})$ and $\eta \in L^0(\mathbb{R}_+^2; \mathbb{F})$, the process $R$ defined by $R_t(x,y) = p_t r(x,y) - \eta_t \cdot (x,y)$ satisfies the above assumptions. Here, $r(x,y)$ stands for the number of units of some good which can be produced by using $x$ (resp. $y$) units of the industrial asset number 1 (resp. 2). $p_t$ denotes the price at time $t$ of this good in terms of the financial asset (interpreted as a currency). $\eta_t$ is the price at time $t$, in terms of the financial asset, of some other goods which are used in the production process. Then, $\eta_t \cdot (x,y)$ stands for a production cost which, here, is assumed to be linear in $(x,y)$.

2.2. *Admissible consumption processes.* A consumption process is an $\mathbb{F}$-adapted process $c = (c_t)_{t \in \mathbb{T}}$ with values in $\mathbb{R}_+^d$. Given an initial endowment $x \in \mathbb{R}^d \times \mathbb{R}_+^N$, we say that a consumption process $c$ is $x$-admissible if $(\sum_{t=0}^T c_t, 0_N) \in A_T(x; K, R)$. We then define

$$(2.2) \quad \mathcal{C}_T(x; K, R) := \left\{ c = (c_t)_{t \in \mathbb{T}} \in L^0(\mathbb{R}_+^d; \mathbb{F}) : \left( \sum_{t \in \mathbb{T}} c_t, 0_N \right) \in A_T(x; K, R) \right\}.$$

Observe that we only allow consumption in terms of financial assets. This formulation is well understood when the financial assets are indeed currencies.

LEMMA 2.1. *Let $(K, R) \in \mathcal{K} \times \mathcal{R}$ be such that (R2) holds and fix $x \in \mathbb{R}^d \times \mathbb{R}_+^N$. Then $A_T(x; K, R)$ is convex, and so is $\mathcal{C}_T(x; K, R)$.*

PROOF. Let $x = (x^F, x^I) \in \mathbb{R}^d \times \mathbb{R}_+^N$, let $g$ and $\tilde{g}$ be two elements of $A_T(x; K, R)$ and let $\xi$ and $\tilde{\xi}$ be two elements of $\mathcal{A}(x; K)$ such that

$$x + \sum_{s=0}^T \xi_s + \sum_{s=0}^{T-1} R_{s+1}(x^I + I(\xi)_s) = g,$$

$$x + \sum_{s=0}^T \tilde{\xi}_s + \sum_{s=0}^{T-1} R_{s+1}(x^I + I(\tilde{\xi})_s) = \tilde{g}.$$

For $\varepsilon \in [0,1]$, we define $\xi^\varepsilon = \varepsilon \xi + (1-\varepsilon)\tilde{\xi}$. Let $\rho^\varepsilon \in L^0(\mathbb{R}^{d+N}; \mathbb{F})$ be defined by $\rho_0^\varepsilon = 0$ and

$$\rho_{t+1}^\varepsilon := \varepsilon R_{t+1}(x^I + I(\xi)_t) + (1-\varepsilon) R_{t+1}(x^I + I(\tilde{\xi})_t) - R_{t+1}(x^I + I(\xi^\varepsilon)_t)$$

for $0 \leq t \leq T-1$. In view of (R2)

$$\rho_t^\varepsilon \in -\underline{K}_t, \qquad t \in \mathbb{T}.$$

Then, $\hat{\xi}^\varepsilon \in L^0(\mathbb{R}^{d+N}; \mathbb{F})$ defined by

$$\hat{\xi}_t^\varepsilon := \varepsilon \xi_t + (1-\varepsilon)\tilde{\xi}_t + \rho_t^\varepsilon, \qquad t \in \mathbb{T},$$

lies in $\mathcal{A}(x; K)$ and satisfies

$$x + \sum_{t=0}^T \hat{\xi}_t^\varepsilon + \sum_{t=0}^{T-1} R_{t+1}(x^I + I(\hat{\xi}^\varepsilon)_t) = \varepsilon g + (1-\varepsilon)\tilde{g}.$$

This concludes the proof. $\square$

**3. The robust no-arbitrage condition.** In order to avoid trivial situations, we need to impose a no-arbitrage condition on the global market. Extending in a natural way the usual notion of no-arbitrage, we assume that

$$NA(K, R): A_T(K, R) \cap L^0(\mathbb{R}_+^{d+N}) = \{0\}.$$

In this paper, we shall indeed impose a stronger condition, which is similar to the one introduced by Schachermayer [13] and further studied by Kabanov, Stricker and Rásonyi [8]. To this end, for $K \in \mathcal{K}$, we define $\underline{K}^0 = (\underline{K}_t^0)_{t \in \mathbb{T}}$ by $\underline{K}_t^0 = \underline{K}_t \cap (-\underline{K}_t)$ for $t \in \mathbb{T}$, and we say that a couple $(\tilde{K}, \tilde{R}) \in \mathcal{K} \times \mathcal{R}$ dominates $(K, R) \in \mathcal{K} \times \mathcal{R}$ if, for each $t \in \mathbb{T}$:

(D1) $\underline{K}_t \setminus \underline{K}_t^0 \subset \mathrm{ri}(\underline{\tilde{K}}_t)$ and $K_t \subset \tilde{K}_t$,
(D2) $\tilde{R}_t(0) \in \underline{K}_t$ and $\tilde{R}_t(\alpha) - R_t(\alpha) \in \mathrm{ri}(\underline{K}_t)$, $\alpha \in \mathbb{R}_+^N \setminus \{0\}$.

We then assume that $(K, R)$ satisfies the *robust no-arbitrage property*:

$NA^r(K, R): NA(\tilde{K}, \tilde{R})$ holds for some $(\tilde{K}, \tilde{R})$ which dominates $(K, R)$.

In the context of pure financial models as in [13] and [8], the robust no-arbitrage condition means that there is no arbitrage even if we slightly reduce the size of the transaction costs which are not already equal to zero. In our context, the same interpretation holds for the financial part of the model. As for the industrial part, we assume that the no-arbitrage property is also stable under a slight increase of the nonlinear returns.

Our first result shows that the $NA^r$ condition implies the closedness of $A_T(x; K, R)$.

THEOREM 3.1. *Let $(K, R) \in \mathcal{K} \times \mathcal{R}$ be such that (R1)–(R3) and $NA^r(K, R)$ hold. Then, for all $x \in \mathbb{R}^d \times \mathbb{R}_+^N$, $A_T(x; K, R)$ and $\mathcal{C}_T(x; K, R)$ are closed in probability.*

PROOF. See Lemma 5.4 and Remark 5.2 below. □

In Section 6 we shall provide a dual formulation for $A_T(x; K, R)$ and $\mathcal{C}_T(x; K, R)$. It is not the main aim of this paper but it will be useful in the proof of Theorem 4.1 below. As usual, the dual formulation is obtained by using the closure property of $A_T(x; K, R)$.

**4. Existence in the utility maximization problem.** We now consider a sequence $(U_t)_{t \in \mathbb{T}}$ of concave mappings from $\mathbb{R}_+^d$ into $\mathbb{R} \cup \{-\infty\}$ such that

$$(4.1) \qquad \mathrm{cl}(\mathrm{dom}(U_t)) = \mathbb{R}_+^d, \qquad t \in \mathbb{T},$$

where $\mathrm{cl}(\mathrm{dom}(U_t))$ denotes the closure of the effective domain of $U_t$, $\mathrm{dom}(U_t) := \{c \in \mathbb{R}^d : |U_t(c)| < \infty\}$. It is natural to assume that $U_t$ is $\mathbb{R}^d$-nondecreasing in the sense that

$$(4.2) \qquad U_t(x) \geq U_t(y) \qquad \text{if } x - y \in \mathbb{R}_+^d, t \in \mathbb{T}.$$

The utility maximization problem is defined as

$$u(x) := \sup_{c \in \mathcal{C}_T^U(x; K, R)} \mathbb{E}\left[\sum_{t \in \mathbb{T}} U_t(c_t)\right], \qquad x \in \mathbb{R}^d \times \mathbb{R}_+^N,$$

where

$$\mathcal{C}_T^U(x; K, R) := \left\{c \in \mathcal{C}_T(x; K, R) : \left(\sum_{t \in \mathbb{T}} U_t(c_t)\right)^- \in L^1(\mathbb{P})\right\}.$$

REMARK 4.1. We claim that $\mathcal{C}_T^U(x; K, R) \neq \varnothing$ whenever $x \in \mathrm{int}(K_0)$. This follows from the following observations.

1. By assumption $\mathbb{R}_+^{d+N} \setminus \{0\} \subset \mathrm{int}(K_0)$. It follows that $(K_0)^* \setminus \{0\} \subset \mathrm{int}(\mathbb{R}_+^{d+N})$. In particular, for $H_1 = \{y \in \mathbb{R}^{d+N} : y^1 = 1\}$, the set $(K_0)^* \cap H_1$ is compact and there is some $\varepsilon > 0$ such that $y^i \geq \varepsilon$ for all $1 \leq i \leq d + N$ and $y \in (K_0)^* \cap H_1$. Also observe that, for $y \in K_0^*$, $y^1 = 0$ implies $y = 0$.
2. Observe now that $x \in \mathrm{int}(K_0)$ if and only if $y \cdot x > 0$ for all $y \in (K_0)^* \cap H_1$. It then follows from part 1 that, for $x \in \mathrm{int}(K_0)$, we can find some $\tilde{x} \in \mathbb{R}^{d+N}$ with $\tilde{x}^i > 0$, for $i \leq d$, such that $x - \tilde{x} \in K_0$.
3. Letting $x$ and $\tilde{x}$ be as in part 2, we define the process $c$ as $c_t^i = \tilde{x}^i / T$ for all $i \leq d$ and $t \in \mathbb{T}$. Then, $c \in \mathcal{C}_T(\tilde{x}; K, R) \subset \mathcal{C}_T(x; K, R)$ and $c_t \in \mathrm{dom}(U_t)$ for all $t \in \mathbb{T}$; see (4.1).

As usual, we need to impose some additional conditions on the utility functions. In our multivariate framework, it is natural to rewrite the usual Inada conditions in terms of the Fenchel transforms associated to $U_t$:

$$\tilde{U}_t(y) = \sup_{x \in \mathbb{R}^d_+} U_t(x) - x \cdot y, \qquad y \in \mathbb{R}^d_+, t \in \mathbb{T}.$$

In the smooth one-dimensional case, the usual Inada conditions $U'_t(0) = +\infty$ and $U'_t(+\infty) = 0$ are equivalent to $\mathrm{dom}(\tilde{U}_t) \supset (0, \infty)$. We therefore assume that

(4.3) $$\mathrm{int}(\mathbb{R}^d_+) \subset \mathrm{dom}(\tilde{U}_t),$$

which is equivalent to $\mathrm{int}(\mathbb{R}^d_+) \subset \bigcup_{x \in (0,\infty)^d} \partial U(x)$, where $\partial U_t(x) := -\partial(-U_t(x))$ and $\partial(-U_t(x))$ denotes the subgradient of $-U_t(x)$ at $x$ in the sense of convex analysis (see, e.g., [12]).

For later use, observe that

(4.4) $$\tilde{U}_t(x) \leq \tilde{U}_t(y) \qquad \text{if } x - y \in \mathbb{R}^d_+, t \in \mathbb{T}.$$

REMARK 4.2. Set $e = (1/\sqrt{d}, \ldots, 1/\sqrt{d}) \in \mathbb{R}^d$ and define the convex map $V_t : r \in \mathbb{R}_+ \mapsto \tilde{U}_t(re)$. By (4.4), $V_t$ is nonincreasing. We claim that

$$\lim_{r \to +\infty} V'_t(r) = 0,$$

where $V'_t$ denotes the right-hand side derivative of $V_t$. Indeed, since $U_t$ and $\tilde{U}_t$ are conjugate functions, it follows from assumption (4.1) that, for all $x \in (0, \infty)^d$,

$$-\infty < U_t(x) = \inf_{y \in (0,\infty)^d} \tilde{U}_t(y) + x \cdot y \leq \inf_{r > 0} V_t(r) + (x \cdot e)r.$$

But if $\lim_{r \to +\infty} V'_t(r) \leq -\varepsilon$ for some $\varepsilon > 0$, then, for $x$ such that $x \cdot e \leq \varepsilon/2$, the right-hand side term is $-\infty$, thus leading to a contradiction.

We finally impose that for all $t \in \mathbb{T}$ and $\lambda \in (0,1)$

(4.5) $$\tilde{U}_t(\lambda y) \leq C^\lambda_t (1 + \tilde{U}_t(y)^+) \qquad \text{for all } y \in \mathrm{dom}(\tilde{U}_t),$$

for some finite $C^\lambda_t \geq 0$.

REMARK 4.3. (i) Assume that $d = 1$, that $U_t$ is strictly concave, increasing, continuously differentiable and that, without loss of generality, $U_t(\infty) > 0$. Then, it follows from Lemma 6.3 in [9] that the *reasonable asymptotic elasticity condition*

(4.6) $$AE(U_t) := \limsup_{x \to \infty} \frac{x U'_t(x)}{U_t(x)} < 1$$

is equivalent to the existence of some $y_0 > 0$ such that, for all $\lambda \in (0,1)$, there is a finite $C_t^\lambda \geq 0$ for which

$$\tilde{U}_t(\lambda y) \leq C_t^\lambda \tilde{U}_t(y) \qquad \text{for all } y < y_0.$$

Since $\tilde{U}_t$ is nonincreasing this implies (4.5).

(ii) In (i), (4.6) can be equivalently written as

$$\limsup_{y \to 0} \frac{-y \tilde{U}_t'(y)}{\tilde{U}_t(y)} < \infty;$$

see Proposition 4.1 in [2]. This condition was extended to the nonsmooth multivariate case by Deelstra, Pham and Touzi [2]. In our setting, their condition takes the form

(4.7) $$\limsup_{\ell(y) \to 0} \left( \sup_{q \in -\partial \tilde{U}_t(y)} q \cdot y \right) \Big/ \tilde{U}_t(y) < \infty,$$

where $\partial \tilde{U}_t(y)$ denotes the subgradient of $\tilde{U}_t$ at $y$ in the sense of convex analysis and

$$\ell(y) := \inf_{x \in \mathbb{R}_+^d, \|x\|=1} x \cdot y.$$

Arguing as in the proof of Lemma 4.1 in [2], one easily checks that (4.7) implies (4.5). However, (4.7) is much stronger than (4.5). For instance, if we take $U_t(x) := \ln(x^1)$, then elementary computations show that (4.5) is satisfied while (4.7) is not.

EXAMPLE 4.1. Let us consider two examples of families of utility functions satisfying the above assumptions. Let $(U_t^i)_{i \leq d}$ be a family of strictly concave, increasing, continuously differentiable utility functions defined on $(0, \infty)$ such that, without loss of generality, $U_t^i(\infty) > 0$ for all $i \leq d$. Assume further that each of them satisfies (4.1), (4.3) and (4.6). As observed in the previous remark, it follows from Lemma 6.3 in [9] that (4.5) holds for each of the corresponding Fenchel transform $\tilde{U}_t^i$.

1. Assume that $U_t^i(0) = 0$ for $i \leq d$, so that $\tilde{U}_t^i \geq 0$, and set

$$U_t(x) := \sum_{i=1}^d U_t^i(x^i) \qquad \text{for } x \in \mathbb{R}_+^d.$$

Then, $\tilde{U}_t(y) = \sum_{i=1}^d \tilde{U}_t^i(y^i)$. Since (4.5) holds for each $\tilde{U}_t^i \geq 0$, it is satisfied by $\tilde{U}_t$ too.

2. Set

$$U_t(x) := U_t^1(x^1) \qquad \text{for } x \in \text{dom}(U_t^1) \times \mathbb{R}_+^{d-1}.$$

Then, $\tilde{U}_t(y) = \tilde{U}_t^1(y^1)$ which also satisfies (4.5).

We can now state our main result.

THEOREM 4.1.  *Fix $(K,R) \in \mathcal{K} \times \mathcal{R}$ such that (R1)–(R3) and $NA^r(K,R)$ hold. Let the conditions* (4.1), (4.2), (4.3) *and* (4.5) *hold. Assume further that $u(x) < \infty$ for some $x \in \mathrm{int}(K_0)$. Then:*

  (i) *$u(x) < \infty$ for all $x \in \mathbb{R}^d \times \mathbb{R}_+^N$,*
  (ii) *for all $x \in \mathbb{R}^d \times \mathbb{R}_+^N$ such that $\mathcal{C}_T^U(x;K,R) \neq \varnothing$, there is some $c^* \in \mathcal{C}_T^U(x;K,R)$ such that*

$$u(x) = \mathbb{E}\left[\sum_{t \in \mathbb{T}} U_t(c_t^*)\right].$$

REMARK 4.4.  If the $U_t$'s are assumed to be strictly concave, then uniqueness holds for the utility maximization problem.

REMARK 4.5.  In Remark 7.3 below, we discuss the assumption (4.5) which can be replaced by a finiteness condition on some auxiliary dual problem as in [10].

REMARK 4.6.  In Kabanov and Kijima [5], the optimization problem is split into an optimal investment problem plus an optimal consumption problem, given the optimal investment strategy. This separation principle is intimately related to the particular "no-bankruptcy" constraint imposed on the wealth process in their complete market framework. In our incomplete market framework and under the natural constraint induced by the definition of the admissibility set $\mathcal{C}_T(x;K,R)$ [see (2.2)], it does not seem possible to obtain such a separation principle.

The remaining sections are organized as follows. In Section 5 we show that $A_T(x;K,R)$ is closed in probability as soon as (R1)–(R3) and $NA^r(K,R)$ hold. In Section 6 we use this result to provide a dual formulation for the set of attainable claims. The proof of Theorem 4.1 is given in Section 7.

**5. The closure property.**  Observe that, because of the constraint (2.1), the sets $A_T(x;K,R)$ are not $K_T$-solid, that is,

$$A_T(x;K,R) \not\supseteq A_T(x;K,R) - L^0(K_T).$$

Indeed, $f \notin A_T(x;K,R)$ whenever $\mathbb{P}[f \notin \mathbb{R}^d \times \mathbb{R}_+^N] > 0$. In order to obtain a suitable dual formulation for $A_T(x;K,R)$ (see Section 6 below), we therefore introduce the $K_T$-solid envelope of $A_T(x;K,R)$:

$$A_T^s(x;K,R) := A_T(x;K,R) - L^0(K_T).$$

Since

$$A^s_T(x; K, R) \cap L^0(\mathbb{R}^d \times \mathbb{R}^N_+) = A_T(x; K, R), \tag{5.1}$$

passing from $A^s_T(x; K, R)$ to $A_T(x; K, R)$ is straightforward. In particular, if $A^s_T(x; K, R)$ is closed in probability, then so is $A_T(x; K, R)$.

In this section, we prove the closedness of $A^s_T(x; K, R)$. It is not of direct use for the proof of Theorem 4.1, that is, the closedness of $A_T(x; K, R)$ is enough, but it will allow us to establish a general dual formulation for the set of elements $g$ of $A_T(x; K, R)$ which are "bounded from below"; see Theorem 6.2 in the next section.

Observe that we can rewrite $A^s_T(x; K, R)$ as

$$A^s_T(x; K, R) = \left\{ x + \sum_{t=0}^T \xi_t + \sum_{t=0}^{T-1} R_{t+1}(x^I + I(\xi)_t), \xi \in \mathcal{A}^s(x; K) \right\},$$

where, for $x = (x^F, x^I) \in \mathbb{R}^d \times \mathbb{R}^N_+$, $\mathcal{A}^s(x; K)$ is the set of adapted processes $\xi$ such that

$$(5.2) \quad \xi_t \in -K_t \quad \text{and} \quad I(\xi)_t \mathbb{1}_{t \leq T-1} + x^I \in \mathbb{R}^N_+, \qquad \mathbb{P}\text{-a.s. for all } t \in \mathbb{T}.$$

We shall simply write $A^s_T(K, R)$ and $\mathcal{A}^s(K)$ when $x = 0$.

The closure property of $A^s_T(x; K, R)$ is a consequence of the stronger result stated in Lemma 5.4 below. We essentially follow the steps of [8]. First, we establish Lemma 5.1 which has to be compared to Lemma 5 in [8]. This is the key result to prove the closure property. We then use an induction argument similar to that in the above paper. The key Lemma 5.1 is used to "transform" unbounded sequences $(\xi^n)$ into bounded ones for which we can extract convergent random subsequences (see Lemma 5.2). Lemma 5.3 is used to control the induced nonlinear returns $R_t$.

LEMMA 5.1. *Let* $(K, R) \in \mathcal{K} \times \mathcal{R}$ *be such that* $NA^r(K, R)$ *holds. Let* $\xi \in \mathcal{A}^s(K)$ *be such that*

$$\sum_{t=0}^T \xi_t + \sum_{t=0}^{T-1} R_{t+1}(I(\xi)_t) = \epsilon$$

*for some* $\epsilon \in \underline{K}_{t_0}$ *with* $t_0 \in \mathbb{T}$. *Then,* $\epsilon \in \underline{K}^0_{t_0}$, *and*

$$I(\xi)_t = 0, \qquad \xi_t \in \underline{K}^0_t \text{ for all } t \in \mathbb{T}.$$

PROOF.  1. First assume that $\mathbb{P}[\epsilon \notin \underline{K}^0_{t_0}] > 0$. By (D1), there is a set $B \subset \Omega$ of positive probability on which $\epsilon \in \text{ri}(\underline{\tilde{K}}_{t_0})$. Hence, we can find some $\beta \in L^0(\mathbb{R}_+^{d+N}; \mathcal{F}_{t_0}) \setminus \{0\}$, such that $-\epsilon + \beta \in -\underline{\tilde{K}}_{t_0}$ on $B$. Set $\hat{\xi}_t = \xi_t + (\beta - \epsilon)\mathbb{1}_{t=t_0}$. Since $\beta - \epsilon$ takes values in $\mathbb{R}^{d+N}$, we have $I(\hat{\xi}) = I(\xi)$ and

$$\sum_{t=0}^{T} \hat{\xi}_t + \sum_{t=0}^{T-1} R_{t+1}(I(\hat{\xi})_t) = \beta.$$

Set

$$r_{t+1} = \tilde{R}_{t+1}(I(\hat{\xi})_t) - R_{t+1}(I(\hat{\xi})_t)$$
$$\tilde{\xi}_0 = \hat{\xi}_0 \quad \text{and} \quad \tilde{\xi}_{t+1} = \hat{\xi}_{t+1} - r_{t+1}, \qquad 0 \leq t \leq T-1.$$

By (D2), $r_{t+1} \in \underline{K}_t$ $\mathbb{P}$-a.s. and $\tilde{\xi} \in \mathcal{A}^s(K) \subset \mathcal{A}^s(\tilde{K})$ satisfies

$$\sum_{t=0}^{T} \tilde{\xi}_t + \sum_{t=0}^{T-1} \tilde{R}_{t+1}(I(\tilde{\xi})_t) = \beta.$$

Since $\beta \in L^0(\mathbb{R}_+^{d+N}; \mathcal{F}_{t_0}) \setminus \{0\}$, this contradicts $NA(\tilde{K}, \tilde{R})$ and therefore $NA^r(K, R)$.

2. If $\mathbb{P}[I(\xi)_{t^*} \neq 0] > 0$ for some $t^* \in \mathbb{T} \setminus \{T\}$, then on a set $B \subset \Omega$ of positive probability we have $I(\xi)_{t^*} \neq 0$. Set $\alpha := \tilde{R}_{t^*+1}(I(\xi)_{t^*}) - R_{t^*+1}(I(\xi)_{t^*})$. Then, by (D2), $\alpha \in \underline{K}_{t^*+1}$ $\mathbb{P}$-a.s. and $\alpha \in \mathrm{ri}(\underline{K}_{t^*+1})$ on $B$. We can then find some $\beta \in L^0(\mathbb{R}_+^{d+N}; \mathcal{F}_{t^*+1}) \setminus \{0\}$ such that $\alpha - \beta \in \underline{K}_{t^*+1}$. Then,

$$-\alpha + \beta + \sum_{t=0}^{T} \xi_t + \sum_{t=0}^{T-1} \tilde{R}_{t+1}(I(\xi)_t) = \epsilon + \beta + \gamma,$$

where

$$\gamma := \sum_{t \in \mathbb{T} \setminus \{t^*\}} \tilde{R}_{t+1}(I(\xi)_t) - R_{t+1}(I(\xi)_t) \in \sum_{t \in \mathbb{T} \setminus \{t^*\}} \underline{K}_t, \qquad \mathbb{P}\text{-a.s.}$$

by (D2). Arguing as in part 1, we obtain a contradiction to $NA(\tilde{K}, \tilde{R})$. Hence, $I(\xi)_t = 0$ $\mathbb{P}$-a.s. for all $t < T$. Since $\epsilon$ takes values in $\mathbb{R}^{d+N}$, we must also have $I(\xi)_T = 0$ $\mathbb{P}$-a.s.

3. We already know from part 2 that $I(\xi)_t = 0$ for each $t \in \mathbb{T}$. It follows that $\xi_t \in -\underline{K}_t$ for all $t \in \mathbb{T}$. Assume that $\mathbb{P}[\xi_{t^*} \notin \underline{K}_{t^*}^0] > 0$ for some $t^* \in \mathbb{T}$. By (D1), there is a set $B \subset \Omega$ of positive probability on which we have $\xi_{t^*} \in -\mathrm{ri}(\underline{\tilde{K}}_{t^*})$. We can then find some $\beta \in L^0(\mathbb{R}_+^{d+N}; \mathcal{F}_{t^*}) \setminus \{0\}$ such that $\xi_{t^*} + \beta \in -\underline{\tilde{K}}_{t^*}$. Since

$$\beta + \sum_{t=0}^{T} \xi_t + \sum_{t=0}^{T-1} R_{t+1}(I(\xi)_t) = \beta + \epsilon,$$

we obtain a contradiction to $NA(\tilde{K}, \tilde{R})$ by the same arguments as in part 1.
□

Before going on with the proof of the closure property, we recall the following lemma whose proof can be found in [6].

LEMMA 5.2. *Set $\mathcal{G} \subset \mathcal{F}$ and let $E$ be a closed subset of $\mathbb{R}^{d+N}$. Let $(\eta^n)_{n \geq 1}$ be a sequence in $L^0(E; \mathcal{G})$. Set $\tilde{\Omega} := \{\liminf_{n \to \infty} \|\eta^n\| < \infty\}$. Then, there is an increasing sequence of random variables $(\tau(n))_{n \geq 1}$ in $L^0(\mathbb{N}; \mathcal{G})$ such that $\tau(n) \to \infty$ $\mathbb{P}$-a.s. and, for each $\omega \in \tilde{\Omega}$, $\eta^{\tau(n)}(\omega)$ converges to some $\eta^*(\omega)$ with $\eta^* \in L^0(E; \mathcal{G})$.*

As a consequence, we first obtain an additional property on $R$ which will be useful in the proof of Lemma 5.4 below.

LEMMA 5.3. *Let $R \in \mathcal{R}$ be such that (R1)–(R3) hold. Let $(\eta^n, \alpha^n)_{n \geq 1}$ be a sequence in $L^0(\mathbb{R}_+ \times \mathbb{R}_+^N; \mathcal{F}_t)$ such that $(\eta^n, \alpha^n) \to (\infty, \alpha)$ $\mathbb{P}$-a.s. for some $\alpha \in L^0(\mathbb{R}_+^N)$. Then, there is a sequence $(\tau_n)_{n \geq 1}$ in $L^0(\mathbb{N}; \mathcal{F}_t)$ such that $\tau_n \to \infty$ $\mathbb{P}$-a.s. and*

$$\lim_{n \to \infty} (\eta^{\tau_n})^{-1} R_t(\eta^{\tau_n} \alpha^{\tau_n}) - R_t(\alpha) = -\epsilon$$

*for some $\epsilon \in L^0(\underline{K}_t; \mathcal{F}_t)$.*

PROOF. By (R1), (R2),

(5.3) $\quad (\eta^n)^{-1} R_t(\eta^n \alpha^n) - R_t(\alpha^n) \in -\underline{K}_t \quad$ on $\{\eta^n \geq 1\}$.

(i) We claim that we can find some $Y \in L^\infty(K_t^*)$ with $Y^i > 0$ $\mathbb{P}$-a.s. for all $i = 1, \ldots, d+N$. Then, on $\{\eta^n \geq 1\}$,

$$Y \cdot [(\eta^n)^{-1} R_t(\eta^n \alpha^n) + (\eta^n)^{-1} a_t + L_t(\alpha^n)]$$
$$\leq Y \cdot [R_t(\alpha^n) + (\eta^n)^{-1} a_t + L_t(\alpha^n)],$$

where $a_t \in L^0(\mathbb{R}_+^{d+N})$ and $L_t \in \mathcal{R}$ are given by (R3). Since $R_t(\alpha^n)$ converges $\mathbb{P}$-a.s. to $R_t(\alpha)$ [see (R1)], $(\eta^n)^{-1} a_t + L_t(\alpha^n)$ converges $\mathbb{P}$-a.s. to $L_t(\alpha)$ and $(\eta^n)^{-1} R_t(\eta^n \alpha^n) + (\eta^n)^{-1} a_t + L_t(\alpha^n) \in L^0(\mathbb{R}_+^{d+N})$, we deduce that

$$\liminf_{n \to \infty} \|(\eta^n)^{-1} R_t(\eta^n \alpha^n)\| < \infty.$$

In view of Lemma 5.2, we can then find a sequence $(\tau_n)_{n \geq 1}$ in $L^0(\mathbb{N}; \mathcal{F}_t)$ such that $\tau_n \to \infty$ $\mathbb{P}$-a.s. and $(\eta^{\tau_n})^{-1} R_t(\eta^{\tau_n} \alpha^{\tau_n})$ converges $\mathbb{P}$-a.s. Since $\underline{K}_t$ is closed, the result then follows from (5.3).

(ii) It remains to prove that we can find some $Y \in L^\infty(K_t^*)$ with $Y^i > 0$ $\mathbb{P}$-a.s. for all $i = 1, \ldots, d+N$. Observe that for $X \in L^0(\mathrm{ri}(K_t))$ there is some $Y \in L^\infty(K_t^*)$ such that $Y \cdot X > 0$. Let $e_i$ be the vector of $\mathbb{R}^{d+N}$ defined by $e_i^j = \mathbb{1}_{i=j}$. Since $K_t$ dominates $\mathbb{R}_+^{d+N}$, that is, $\mathbb{R}_+^{d+N} \setminus \{0\} \subset \mathrm{ri}(K_t)$, for each $1 \leq i \leq d+N$ we can find some $Y_i \in K_t^*$ such that $Y_i \cdot e_i > 0$. Then, $Y := \sum_{i=1}^{d+N} Y_i \in K_t^*$ satisfies the required property. $\square$

REMARK 5.1. In the above proof, assumption (R3) was used only to show that

(5.4) $$\liminf_{n\to\infty} \|(\eta^n)^{-1} R_t(\eta^n \alpha^n)\| < \infty, \qquad \mathbb{P}\text{-a.s.}$$

Then, we could replace (R3) by: for all sequence $(\eta^n, \alpha^n)_{n\geq 1}$ in $L^0(\mathbb{R}_+ \times \mathbb{R}_+^N; \mathcal{F}_t)$ such that $(\eta^n, \alpha^n) \to (\infty, \alpha)$ $\mathbb{P}$-a.s. for some $\alpha \in L^0(\mathbb{R}_+^N)$, we have (5.4).

We can now state the main result of this section.

LEMMA 5.4. *Let* $(K, R) \in \mathcal{K} \times \mathcal{R}$ *be such that* (R1)–(R3) *and* $NA^r(K, R)$ *hold. For* $t \in \mathbb{T}$ *and* $\alpha \in L^0(\mathbb{R}_+^N; \mathcal{F}_t)$, *let* $\mathcal{Y}^{t,\alpha}(K)$ *be the set of processes* $\xi \in L^0(\mathbb{R}^{d+N}; \mathbb{F})$ *such that*

$$\xi_s \in -K_s \mathbb{1}_{s\geq t} \quad \text{for all } s \in \mathbb{T} \quad \text{and}$$
$$I(\xi)_s + \alpha \in \mathbb{R}_+^N \quad \text{for all } t \leq s \leq T-1, \ \mathbb{P}\text{-a.s.}$$

*For* $t \in \mathbb{T}$, *let* $Y_T^t(K, R)$ *denote the set of elements* $(\alpha, g) \in L^0(\mathbb{R}_+^N; \mathcal{F}_t) \times L^0(\mathbb{R}^{d+N}; \mathcal{F}_T)$ *such that there is some* $\xi \in \mathcal{Y}^{t,\alpha}(K)$ *for which*

$$\sum_{s=t}^{T} \xi_s + \sum_{s=t}^{T-1} R_{s+1}(I(\xi)_s + \alpha) = g.$$

*Then, for all* $t \in \mathbb{T}$, $Y_T^t(K, R)$ *is closed for the convergence in probability.*

REMARK 5.2. For $x = (x^F, x^I) \in \mathbb{R}^d \times \mathbb{R}_+^N$, the above lemma readily implies that $A_T^s(x; K, R)$ is closed in probability since $(x^I, g_n + (x^F, 0_N)) \in Y_T^0(K, R)$ if and only if $g_n \in A_T^s(x; K, R)$. In view of (5.1), this shows that $A_T(x; K, R)$ is closed too and so is $\mathcal{C}_T(x; K, R)$.

PROOF. We proceed by induction. For $t = T$, there is nothing to prove. We then assume that $Y_T^{t+1}(K, R)$ is closed for some $0 \leq t < T$ and show that this implies that $Y_T^t(K, R)$ is closed also. Let $(\alpha^n, g^n)_{n\geq 1}$ be a sequence in $Y_T^t(K, R)$ that converges in probability to some $(\alpha, g) \in L^0(\mathbb{R}_+^N; \mathcal{F}_t) \times L^0(\mathbb{R}^{d+N}; \mathcal{F}_T)$. After passing to a subsequence, we can assume that the convergence holds $\mathbb{P}$-a.s. Let $(\xi^n)_{n\geq 1}$ be a sequence such that

(5.5) $\xi^n \in \mathcal{Y}^{t,\alpha^n}(K)$ and $\sum_{s=t}^{T} \xi_s^n + \sum_{s=t}^{T-1} R_{s+1}(I(\xi^n)_s + \alpha^n) = g^n, \qquad n \geq 1.$

Set $\tilde{\Omega} = \{\liminf_{n\to\infty} \|\xi_t^n\| < \infty\}$ and observe that $\tilde{\Omega} \in \mathcal{F}_t$.

1. By Lemma 5.2, if $\mathbb{P}[\tilde{\Omega}] = 1$, we can find an increasing sequence of random variables $(\tau(n))_{n\geq 1}$ in $L^0(\mathbb{N}; \mathcal{F}_t)$ such that, for each $\omega \in \tilde{\Omega}$, $\xi_t^{\tau(n)}(\omega)$ converges to some $\xi_t(\omega)$ with $\xi_t \in L^0(\mathbb{R}^{d+N};\mathcal{F}_t)$. We then have

$$g^{\tau(n)} = \xi_t^{\tau(n)} + R_{t+1}(I(\xi^{\tau(n)})_t + \alpha^{\tau(n)})$$
$$+ \sum_{s=t+1}^{T} \tilde{\xi}_s^{\tau(n)} + \sum_{s=t+1}^{T-1} R_{s+1}(I(\tilde{\xi}^{\tau(n)})_s + I(\xi^{\tau(n)})_t + \alpha^{\tau(n)}),$$

where

$$\tilde{\xi}_s^{\tau(n)} = \xi_s^{\tau(n)} \mathbb{1}_{s \geq t+1}, \qquad 0 \leq s \leq T.$$

Hence, $(I(\xi^{\tau(n)})_t + \alpha^{\tau(n)}, g^{\tau(n)} - \xi_t^{\tau(n)} - R_{t+1}(I(\xi^{\tau(n)})_t + \alpha^{\tau(n)}))$ belongs to $Y_T^{t+1}(K,R)$. Since $Y_T^{t+1}(K,R)$ is closed, we can find some $\tilde{\xi} \in L^0(\mathbb{R}^{d+N};\mathbb{F})$, with $\tilde{\xi}_s = 0$ for $s < t+1$, such that

$$\sum_{s=t+1}^{T} \tilde{\xi}_s + \sum_{s=t+1}^{T-1} R_{s+1}(I(\tilde{\xi})_s + I(\xi)_t + \alpha) = g - \xi_t - R_{t+1}(I(\xi)_t + \alpha),$$

where we used (R1) to pass to the limit in $R_{t+1}$. Set

$$\bar{\xi}_s = \xi_t \mathbb{1}_{\{s=t\}} + \tilde{\xi} \mathbb{1}_{\{t < s \leq T\}}, \qquad s \in \mathbb{T}.$$

Then,

$$\sum_{s=t}^{T} \bar{\xi}_s + \sum_{s=t}^{T-1} R_{s+1}(I(\bar{\xi})_s + \alpha) = g$$

where, in view of (5.5),

$$\bar{\xi}_s \in -K_s \mathbb{1}_{s \geq t} \qquad \text{for } s \in \mathbb{T} \quad \text{and}$$
$$I(\bar{\xi})_s + \alpha \in \mathbb{R}_+^N \qquad \text{for } t \leq s \leq T-1, \ \mathbb{P}\text{-a.s.}$$

This shows that $(\alpha, g) \in Y_T^t(K,R)$.

2. We next consider the case where $\mathbb{P}[\tilde{\Omega}] < 1$. Since $\tilde{\Omega} \in \mathcal{F}_t$, we can work separately on $\tilde{\Omega}$ and $\tilde{\Omega}^c$, by considering two alternative strategies depending on the occurrence of $\tilde{\Omega}$ or $\tilde{\Omega}^c$. We can then proceed as if $\mathbb{P}[\tilde{\Omega}^c] = 1$.

2a. Let $\eta_t^n := \|\xi_t^n\| + 1$. Since $\liminf_{n \to \infty} (\eta_t^n)^{-1} \|\xi_t^n\| < \infty$ $\mathbb{P}$-a.s., we can find an increasing sequence of random variables $(\tau(n))_{n \geq 1}$ in $L^0(\mathbb{N};\mathcal{F}_t)$ such that for each $\omega \in \tilde{\Omega}^c$, $(\eta_t^{\tau(n)})^{-1} \xi_t^{\tau(n)}$ converges to some $\bar{\xi}_t^*$ in $L^0(\mathbb{R}^{d+N};\mathcal{F}_t)$.

Set

$$(\bar{\xi}^n, \bar{g}^n, \bar{\alpha}^n) := (\eta_t^{\tau(n)})^{-1}(\xi^{\tau(n)}, g^{\tau(n)}, \alpha^{\tau(n)}) \quad \text{and} \quad \bar{\eta}_t^n := \eta_t^{\tau(n)},$$

so that

(5.6)
$$\bar{g}^n = \bar{\xi}^n_t + (\bar{\eta}^n_t)^{-1} R_{t+1}(\bar{\eta}^n_t(I(\bar{\xi}^n)_t + \bar{\alpha}^n))$$
$$+ \sum_{s=t+1}^{T} \bar{\xi}^n_s + \sum_{s=t+1}^{T-1} (\bar{\eta}^n_t)^{-1} R_{s+1}(\bar{\eta}^n_t(I(\bar{\xi}^n)_s + \bar{\alpha}^n)).$$

Set

(5.7)
$$r^n_{s+1} := R_{s+1}(I(\bar{\xi}^n)_s + \bar{\alpha}^n) - (\bar{\eta}^n_t)^{-1} R_{s+1}(\bar{\eta}^n_t(I(\bar{\xi}^n)_s + \bar{\alpha}^n)),$$
$$t+1 \leq s \leq T-1.$$

In view of (R1), (R2), $r^n_{s+1} \in \underline{K}_{s+1}$, $t+1 \leq s \leq T-1$, $\mathbb{P}$-a.s. Set

(5.8) $\quad \tilde{\xi}^n_s := \bar{\xi}^n_s \mathbb{1}_{s \geq t+1} - r^n_s \mathbb{1}_{s \geq t+2} \in -K_s, \quad s \in \mathbb{T}.$

Since $I(\xi)$ does not depend on the $d$ first component of $\xi$, we have

(5.9) $\quad I(\bar{\xi}^n)_s = I(\tilde{\xi}^n)_s + I(\bar{\xi}^n)_t, \quad s \geq t+1.$

Since $\bar{\alpha}^n + I(\bar{\xi}^n)_t \to I(\bar{\xi}^*)_t$ $\mathbb{P}$-a.s., we deduce from Lemma 5.3 that there is some $\epsilon \in L^0(\underline{K}_{t+1}; \mathcal{F}_{t+1})$ and an increasing sequence of random variables $(\sigma(n))_{n \geq 1}$ in $L^0(\mathbb{N}; \mathcal{F}_{t+1})$ such that

(5.10) $\quad \lim_{n \to \infty} (\bar{\eta}^{\sigma(n)}_t)^{-1} R_{t+1}(\bar{\eta}^{\sigma(n)}_t(I(\bar{\xi}^{\sigma(n)})_t + \bar{\alpha}^{\sigma(n)})) - R_{t+1}(I(\bar{\xi}^*)_t) = -\epsilon,$

where $\sigma(n)$ goes to $\infty$ $\mathbb{P}$-a.s. Since by (5.6)–(5.9)

$$\bar{g}^{\sigma(n)} = \bar{\xi}^{\sigma(n)}_t + (\bar{\eta}^{\sigma(n)}_t)^{-1} R_{t+1}(\bar{\eta}^{\sigma(n)}_t(I(\bar{\xi}^{\sigma(n)})_t + \bar{\alpha}^{\sigma(n)}))$$
$$+ \sum_{s=t+1}^{T} \tilde{\xi}^{\sigma(n)}_s + \sum_{s=t+1}^{T-1} R_{s+1}(I(\tilde{\xi}^{\sigma(n)})_s + I(\bar{\xi}^{\sigma(n)})_t + \bar{\alpha}^{\sigma(n)}),$$

hence, $(I(\bar{\xi}^{\sigma(n)})_t + \bar{\alpha}^{\sigma(n)}, \bar{g}^{\sigma(n)} - (\bar{\eta}^{\sigma(n)}_t)^{-1} R_{t+1}(\bar{\eta}^{\sigma(n)}_t(I(\bar{\xi}^{\sigma(n)})_t + \bar{\alpha}^{\sigma(n)})) - \bar{\xi}^{\sigma(n)}_t)$ belongs to $Y^{t+1}_T(K, R)$. Since $Y^{t+1}_T(K, R)$ is closed and $(\bar{g}^n, \bar{\alpha}^n)$ goes to 0 $\mathbb{P}$-a.s., we can find some adapted process $\tilde{\xi}^*$ such that

$$\tilde{\xi}^*_s \in -K_s \mathbb{1}_{s \geq t+1} \quad \text{for } s \in \mathbb{T},$$
$$I(\tilde{\xi}^*)_s + I(\bar{\xi}^*)_t \in \mathbb{R}^N_+ \quad \text{for all } s \in \mathbb{T} \setminus \{T\}, \ \mathbb{P}\text{-a.s.},$$
$$0 = \lim_{n \to \infty} (\bar{\eta}^{\sigma(n)}_t)^{-1} R_{t+1}(\bar{\eta}^{\sigma(n)}_t(I(\bar{\xi}^{\sigma(n)})_t + \bar{\alpha}^{\sigma(n)})) + \bar{\xi}^*_t + \sum_{s=t+1}^{T} \tilde{\xi}^*_s$$
$$+ \sum_{s=t+1}^{T-1} R_{s+1}(I(\tilde{\xi}^*)_s + I(\bar{\xi}^*)_t),$$

and it follows from (5.10) that

$$L^0(\underline{K}_{t+1}; \mathcal{F}_{t+1}) \ni \epsilon = R_{t+1}(I(\bar{\xi}^*)_t)$$
$$+ \sum_{s=t+1}^{T} \tilde{\xi}_s^* + \bar{\xi}_t^* + \sum_{s=t+1}^{T-1} R_{s+1}(I(\tilde{\xi}^*)_s + I(\bar{\xi}^*)_t).$$

We then define

(5.11) $$\hat{\xi}_s^* := \bar{\xi}_t^* \mathbb{1}_{s=t} + \tilde{\xi}_s^* \mathbb{1}_{s \geq t+1}, \qquad s \in \mathbb{T}.$$

With this new notation, we have $I(\hat{\xi}^*)_s \mathbb{1}_{s \leq T-1} \in \mathbb{R}_+^N$, $\hat{\xi}_s^* \in -K_s \mathbb{1}_{s \geq t}$ for all $s \in \mathbb{T}$, and

(5.12) $$\epsilon = \sum_{s=t}^{T} \hat{\xi}_s^* + \sum_{s=t}^{T-1} R_{s+1}(I(\hat{\xi}^*)_s).$$

By Lemma 5.1, we must have $\epsilon \in \underline{K}_{t+1}^0$,

(5.13) $$I(\hat{\xi}^*)_s = 0 \quad \text{and} \quad \hat{\xi}_s^* \in \underline{K}_s^0 \qquad \text{for all } s \in \mathbb{T}.$$

Finally, letting

$$\check{\xi}_s^* := \hat{\xi}_s^* - \epsilon \mathbb{1}_{s=t+1}, \qquad s \in \mathbb{T},$$

we deduce from (5.11)–(5.13) and (R1) that

(5.14) $$\check{\xi}^* \in \mathcal{A}^s(K), \qquad \check{\xi}_s^* \in -\underline{K}_s \mathbb{1}_{s \geq t} \qquad \text{for all } s \in \mathbb{T} \text{ and } \sum_{s=t}^{T} \check{\xi}_s^* = 0.$$

2b. Since $\|\bar{\xi}_t^*\| = \|\check{\xi}_t^*\| = 1$ on $\tilde{\Omega}$, there is a partition of $\tilde{\Omega}$ into disjoint subsets $\Gamma_i \in \mathcal{F}_t$ such that $\Gamma_i \subset \{(\check{\xi}_t^*)^i \neq 0\}$ for $i = 1, \ldots, d$. We then define

$$\check{\xi}_s^n = \sum_{i=1}^{d} (\xi_s^n - \beta_t^{n,i} \check{\xi}_s^*) \mathbb{1}_{\Gamma_i}, \qquad s \in \mathbb{T},$$

with $\beta_t^{n,i} = (\xi_t^n)^i / (\check{\xi}_t^*)^i$ on $\Gamma_i$, $i = 1, \ldots, d$. Since, by (5.14) and the definition of $\xi^n$,

$$\sum_{s=t}^{T} \check{\xi}_s^n = \sum_{s=t}^{T} \xi_s^n, \qquad \check{\xi}_s^n \in -K_s \mathbb{1}_{s \geq t} \quad \text{and} \quad I(\check{\xi}^n)_s = I(\xi^n)_s, \qquad s \in \mathbb{T},$$

it follows that $\check{\xi}^n \in \mathcal{Y}^{t,\alpha^n}(K)$ and

$$\sum_{s=t}^{T} \check{\xi}_s^n + \sum_{s=t}^{T-1} R_{s+1}(I(\check{\xi}^n)_s + \alpha^n) = g^n, \qquad n \geq 1.$$

We can then proceed as in [8] and obtain the required result by repeating the above argument with $(\check{\xi}^n)_{n \geq 1}$ instead of $(\xi^n)_{n \geq 1}$ and by iterating this procedure a finite number of times. □

**6. Dual formulation for attainable terminal wealth.** In this section we provide a dual characterization of the set of attainable terminal wealth. To this end, given $K \in \mathcal{K}$ and $\tilde{\mathbb{P}} \sim \mathbb{P}$, we define $\mathcal{Z}_T(K, \tilde{\mathbb{P}})$ as the set of adapted processes $Z = (Z^F, Z^I) \in L^1(\mathbb{R}^{d+N}; \tilde{\mathbb{P}}, \mathbb{F})$ such that:

(i) $(Z_t^F, 0_N) \in \text{ri}((\underline{K}_t)^*)$ for each $t \in \mathbb{T}$ and $Z_T \in (K_T)^* \setminus \{0\}$ $\mathbb{P}$-a.s.,
(ii) $Z^F$ is a $\tilde{\mathbb{P}}$-martingale.

REMARK 6.1. Recall that, by assumption, $\mathbb{R}_+^{d+N} \setminus \{0\} \subset \text{int}(K_T)$ $\mathbb{P}$-a.s. It follows that $(K_T)^* \setminus \{0\} \subset \text{int}(\mathbb{R}_+^{d+N})$ $\mathbb{P}$-a.s. This shows that $Z_T \in \text{int}(\mathbb{R}_+^{d+N})$ $\mathbb{P}$-a.s. whenever $Z \in \mathcal{Z}_T(K, \tilde{\mathbb{P}})$ for some $\tilde{\mathbb{P}} \sim \mathbb{P}$.

We start with a series of lemmas which are similar to results in [13] and [8].

LEMMA 6.1. *Fix $(K, R) \in \mathcal{K} \times \mathcal{R}$ satisfying* (R1)–(R3) *and $NA^r(K, R)$. Then, for all $\tilde{\mathbb{P}} \sim \mathbb{P}$, there is a process $Z \in \mathcal{Z}_T(K, \tilde{\mathbb{P}}) \cap L^\infty$ such that*

$$\sup_{g \in A_T^s(K,R) \cap L^1(\tilde{\mathbb{P}})} \mathbb{E}^{\tilde{\mathbb{P}}}[Z_T \cdot g] < \infty.$$

PROOF. Since, by Remark 5.2, $A_T^s(K, R)$ is closed in probability, $A_T^s(K, R) \cap L^1(\tilde{\mathbb{P}})$ is closed in $L^1(\tilde{\mathbb{P}})$. By Lemma 2.1 it is also convex. In view of $NA(K, R)$, which is trivially implied by $NA^r(K, R)$, it then follows from the Hahn–Banach separation theorem that, for each $\phi \in L^1(\mathbb{R}_+^{d+N}; \tilde{\mathbb{P}}) \setminus \{0\}$, we can find $\eta(\phi) \in L^\infty(\mathbb{R}^{d+N})$ such that

$$\mathbb{E}^{\tilde{\mathbb{P}}}[\eta(\phi) \cdot g] < \mathbb{E}^{\tilde{\mathbb{P}}}[\eta(\phi) \cdot \phi] \qquad \text{for all } g \in A_T^s(K, R) \cap L^1(\tilde{\mathbb{P}}).$$

Since $-L^0(K_T) \subset A_T^s(K, R)$, we must have $\eta(\phi) \in L^0((K_T)^*)$. Using a standard exhaustion argument, we obtain some $\eta \in L^0((K_T)^*)$ such that $\mathbb{E}^{\tilde{\mathbb{P}}}[\eta \cdot g] \leq 0$ for all $g \in A_T^s(K, R) \cap L^1(\tilde{\mathbb{P}})$, and $\mathbb{P}[\eta = 0] = 0$. Set $Z_t = (Z_t^F, Z_t^I) = \mathbb{E}[\eta | \mathcal{F}_t]$. Then, $Z^F$ is a martingale. Since $\sum_{t \in \mathbb{T}} -L^0(\underline{K}_t; \mathcal{F}_t) \subset A_T^s(K, R)$, we must have

$$\mathbb{E}^{\tilde{\mathbb{P}}}[\eta \cdot g] \leq 0 \qquad \text{for all } g \in \sum_{t \in \mathbb{T}} -L^1(\underline{K}_t; \tilde{\mathbb{P}}, \mathcal{F}_t).$$

In particular, this shows that $(Z_t^F, 0_N) \in L^0(\text{ri}((\underline{K}_t)^*))$. The rest of the proof then goes as in Corollary 1 in [8] by using Lemma 5.1 and the fact that the $\underline{K}_t = K_t \cap \underline{\mathbb{R}}^{d+N}$ are countably generated (see the remark after Corollary 1 in [8]). □

REMARK 6.2. Observe that $x \in K_0$ if and only if $y \cdot x \geq 0$ for all $y \in (K_0)^* \cap H_1$, where $H_1 = \{y \in \mathbb{R}^{d+N} : y^1 = 1\}$. Using part 1 of Remark 4.1, we then deduce that, for any $x \in \mathbb{R}^d \times \mathbb{R}_+^N$, we can find some $\hat{x} = (\hat{x}^1, 0_{d-1+N}) \in \mathbb{R}^d \times \mathbb{R}_+^N$ such that $\hat{x} - x \in K_0$.

COROLLARY 6.1. *Fix $(K, R) \in \mathcal{K} \times \mathcal{R}$ such that (R1)–(R3) and $NA^r(K, R)$ hold. Fix $x = (0_d, x^I) \in \mathbb{R}^d \times \mathbb{R}_+^N$. Then, for all $\tilde{\mathbb{P}} \sim \mathbb{P}$, there is some $Z \in \mathcal{Z}_T(K, \tilde{\mathbb{P}}) \cap L^\infty$ such that*

$$a(x^I; Z, \tilde{\mathbb{P}}) := \sup_{g \in A_T^s(x; K, R) \cap L^1(\tilde{\mathbb{P}})} \mathbb{E}^{\tilde{\mathbb{P}}}[Z_T \cdot g] < \infty.$$

PROOF. In view of Remark 6.2, there is some $\hat{x} \in \mathbb{R}^{d+N}$ such that $\hat{x} - x \in K_0$. It follows that $A_T^s(x; K, R) \subset A_T^s(\hat{x}; K, R)$. Then, the required result is a direct consequence of Lemma 6.1. Indeed, we can find some $Z$ which satisfies the assertions of Lemma 6.1. Since $A_T^s(x; K, R) - \hat{x} \subset A_T^s(\hat{x}; K, R) - \hat{x} = A_T^s(K, R)$ (see Remark 2.1), it follows that

$$\sup_{g \in A_T^s(x; K, R) \cap L^1(\tilde{\mathbb{P}})} \mathbb{E}^{\tilde{\mathbb{P}}}[Z_T \cdot (g - \hat{x})] \leq \sup_{g \in A_T^s(K, R) \cap L^1(\tilde{\mathbb{P}})} \mathbb{E}^{\tilde{\mathbb{P}}}[Z_T \cdot g],$$

where $Z_T \cdot \hat{x} \in L^\infty$ since $Z_T \in L^\infty$.  □

LEMMA 6.2. *Fix $(K, R) \in \mathcal{K} \times \mathcal{R}$ such that (R1)–(R3) and $NA^r(K, R)$ hold. Fix $x = (x^F, x^I) \in \mathbb{R}^d \times \mathbb{R}_+^N$, $g \in L^0(\mathbb{R}^{d+N}; \mathcal{F}_T)$ and $\tilde{\mathbb{P}} \sim \mathbb{P}$ such that $g \in L^1(\mathbb{R}^{d+N}; \tilde{\mathbb{P}})$. Then,*

$$\mathbb{E}^{\tilde{\mathbb{P}}}[Z_T \cdot g - Z_0^F \cdot x^F] - a(x^I; Z, \tilde{\mathbb{P}}) \leq 0 \quad \text{for all } Z = (Z^F, Z^I) \in \mathcal{Z}_T(K, \tilde{\mathbb{P}})$$

*implies $g \in A_T^s(x; K, R)$.*

PROOF. Fix some $\tilde{\mathbb{P}}$ such that $g \in L^1(\mathbb{R}^{d+N}; \tilde{\mathbb{P}})$. Assume that $g \notin A_T^s(x; K, R) \cap L^1(\tilde{\mathbb{P}})$. Since, by Lemma 2.1 and Remark 5.2, $A_T^s(x; K, R) \cap L^1(\tilde{\mathbb{P}})$ is closed in $L^1(\tilde{\mathbb{P}})$ and convex, we can find some $\eta \in L^\infty(\mathbb{R}^{d+N})$ such that

$$\sup_{g \in A_T^s(x; K, R) \cap L^1(\tilde{\mathbb{P}})} \mathbb{E}^{\tilde{\mathbb{P}}}[\eta \cdot (g - (x^F, 0_N))] < \mathbb{E}^{\tilde{\mathbb{P}}}[\eta \cdot (g - (x^F, 0_N))].$$

Set $Z_t := \mathbb{E}^{\tilde{\mathbb{P}}}[\eta | \mathcal{F}_t]$. The same argument as in Lemma 6.1 shows that $Z^F$ is a $\tilde{\mathbb{P}}$-martingale with $Z_T \in L^0(K_T^*)$ and $(Z^F, 0_N)_t \in L^0((\underline{K}_t)^*; \mathcal{F}_t)$ for each $t \in \mathbb{T}$. Fix $\hat{Z} \in \mathcal{Z}_T(K, \tilde{\mathbb{P}}) \cap L^\infty$ such that $a(x^I; \hat{Z}, \tilde{\mathbb{P}}) < \infty$ (which is possible by Corollary 6.1). For $\varepsilon > 0$ sufficiently small, we have $Z^\varepsilon := \varepsilon \hat{Z} + (1 - \varepsilon) Z \in \mathcal{Z}_T(K, \tilde{\mathbb{P}})$ and

$$a(x^I; Z^\varepsilon, \tilde{\mathbb{P}}) = \sup_{g \in A_T^s(x; K, R) \cap L^1(\tilde{\mathbb{P}})} \mathbb{E}^{\tilde{\mathbb{P}}}[Z_T^\varepsilon \cdot (g - (x^F, 0_N))]$$

$$< \mathbb{E}^{\tilde{\mathbb{P}}}[Z_T^\varepsilon \cdot (g - (x^F, 0_N))],$$

where we used the fact $A_T^s(x; K, R) - (x^F, 0_N) = A_T^s((0_d, x^I); K, R)$. This leads to a contradiction since $(Z^\varepsilon)^F$ is a martingale.  □

We can now state a first version of the so-called *super-hedging theorem*.

THEOREM 6.1. *Fix $(K, R) \in \mathcal{K} \times \mathcal{R}$ such that (R1)–(R3) and $NA^r(K, R)$ hold. Fix $x = (x^F, x^I) \in \mathbb{R}^d \times \mathbb{R}_+^N$. Then, we have the equivalence between:*

(i) $g \in A_T^s(x; K, R)$,

(ii) *for some $\tilde{\mathbb{P}} \sim \mathbb{P}$ such that $g \in L^1(\mathbb{R}^{d+N}; \tilde{\mathbb{P}})$, we have for each $Z = (Z^F, Z^I) \in \mathcal{Z}_T(K, \tilde{\mathbb{P}})$*

$$\mathbb{E}^{\tilde{\mathbb{P}}}[Z_T \cdot g - Z_0^F \cdot x^F] - a(x^I; Z, \tilde{\mathbb{P}}) \leq 0,$$

(iii) *for all $\tilde{\mathbb{P}} \sim \mathbb{P}$ such that $g \in L^1(\mathbb{R}^{d+N}; \tilde{\mathbb{P}})$, we have for each $Z = (Z^F, Z^I) \in \mathcal{Z}_T(K, \tilde{\mathbb{P}})$*

$$\mathbb{E}^{\tilde{\mathbb{P}}}[Z_T \cdot g - Z_0^F \cdot x^F] - a(x^I; Z, \tilde{\mathbb{P}}) \leq 0.$$

PROOF. Since $Z^F$ is a martingale, (i) implies (iii) by definition of $a(x^I; Z, \tilde{\mathbb{P}})$ and the fact that $(g - (x^F, 0_N)) \in A_T^s((0_d, x^I); K, R)$. Obviously (iii) implies (ii). The implication (ii) $\Rightarrow$ (i) follows from Lemma 6.2. $\square$

In the case where the claim is uniformly bounded from below for the natural partial order induced by $K_T$, we can obtain a version of the *superhedging theorem* which does not depend on the integrability properties of $g$.

THEOREM 6.2. *Fix $(K, R) \in \mathcal{K} \times \mathcal{R}$ such that (R1)–(R3) and $NA^r(K, R)$ hold. Fix $x = (x^F, x^I) \in \mathbb{R}^d \times \mathbb{R}_+^N$ and let $g$ be an element of $L^0(\mathbb{R}^{d+N})$ such that $g + c \in K_T$ for some constant $c \in \mathbb{R}^{d+N}$. Then, we have the equivalence between:*

(i) $g \in A_T^s(x; K, R)$,

(ii) *for each $\tilde{\mathbb{P}} \sim \mathbb{P}$ and $Z = (Z^F, Z^I) \in \mathcal{Z}_T(K, \tilde{\mathbb{P}})$*

$$\mathbb{E}^{\tilde{\mathbb{P}}}[Z_T \cdot g - Z_0^F \cdot x^F] - a(x^I; Z, \tilde{\mathbb{P}}) \leq 0,$$

(iii) *for some $\tilde{\mathbb{P}} \sim \mathbb{P}$, we have for each $Z = (Z^F, Z^I) \in \mathcal{Z}_T(K, \tilde{\mathbb{P}})$*

$$\mathbb{E}^{\tilde{\mathbb{P}}}[Z_T \cdot g - Z_0^F \cdot x^F] - a(x^I; Z, \tilde{\mathbb{P}}) \leq 0.$$

PROOF. 1. Let $g \in A_T^s(x; K, R)$ be such that $g + c \in K_T$ for some constant $c = (c^F, 0_N) \in \mathbb{R}^{d+N}$. For $k \geq 1$, set $B_k = \{\|g + c\| \leq k\}$. Then, $\mathbb{1}_{B_k}$ goes to 1 $\mathbb{P}$-a.s. as $k \to \infty$. For each $k \geq 1$, define $g_k := (g + c)\mathbb{1}_{B_k}$. Since $g + c \in L^0(K_T)$, $g_k \in A_T^s(x + c; K, R) = c^F + A_T^s(x; K, R)$ for all $k \geq 1$. Since $g_k$ is bounded, we deduce from Theorem 6.1 that, for each $\tilde{\mathbb{P}} \sim \mathbb{P}$ and $Z = (Z^F, Z^I) \in \mathcal{Z}_T(K, \tilde{\mathbb{P}})$, we must have

$$\mathbb{E}^{\tilde{\mathbb{P}}}[Z_T \cdot g_k - Z_0^F \cdot (x^F + c^F)] - a(x^I; Z, \tilde{\mathbb{P}}) \leq 0.$$

Since $g_k \in L^0(K_T)$ and $Z_T \in L^0((K_T)^*)$, we have $Z_T \cdot g_k \geq 0$ $\mathbb{P}$-a.s. Using Fatou's lemma, we then deduce that

$$\mathbb{E}^{\tilde{\mathbb{P}}}[Z_T \cdot (g+c) - Z_0^F \cdot (x^F + c^F)] - a(x^I; Z, \tilde{\mathbb{P}}) \leq 0 \quad \text{for all } k \geq 1,$$

and (ii) follows from the martingale property of $Z^F$.

2. To see that (ii) implies (i), we define $\tilde{\mathbb{P}} \sim \mathbb{P}$ by $\tilde{\mathbb{P}} = (e^{-\|g\|}/\mathbb{E}[e^{-\|g\|}]) \cdot \mathbb{P}$. Then, $g \in L^1(\mathbb{R}^{d+N}; \tilde{\mathbb{P}})$ and the result follows from Theorem 6.1.

3. Obviously (ii) implies (iii). It remains to check the converse implication. Fix $\tilde{\mathbb{P}}$ such that (iii) holds, $\hat{\mathbb{P}} \sim \mathbb{P}$ and let $H_t := \mathbb{E}^{\tilde{\mathbb{P}}}[d\hat{\mathbb{P}}/d\tilde{\mathbb{P}}|\mathcal{F}_t]$. Then, for $\hat{Z} \in \mathcal{Z}_T(K, \hat{\mathbb{P}})$, we have $\tilde{Z} := (H_t \hat{Z}_t)_{t \in \mathbb{T}} \in \mathcal{Z}_T(K, \tilde{\mathbb{P}})$ and $a(x^I; \hat{Z}, \hat{\mathbb{P}}) = a(x^I; \tilde{Z}, \tilde{\mathbb{P}})$. This shows that (iii) implies (ii). □

REMARK 6.3. Observe from (5.1) that Theorems 6.1 and 6.2 actually provide a dual formulation for $A_T(x; K, R)$. It suffices to add the condition $g \in L^0(\mathbb{R}^d \times \mathbb{R}^N_+)$.

REMARK 6.4. It is clear from the proofs that the results of this section still hold if we replace (R1)–(R3) by the assumption that $A_T(x; K, R)$ is closed.

REMARK 6.5. Although the dual formulation we obtained is already much more general than what we need for the proof of Theorem 4.1, we think that a more precise description of the natural set of dual variables could be obtained by means of Lemma 5.4, which is actually much stronger than the version we used in the proofs. We leave this point for future research.

**7. Proof of the existence result for the optimal consumption problem.** As already explained in the Introduction, the one-dimensional argument of Kramkov and Schachermayer [10] does not work directly in our multivariate setting. We therefore manage this difficulty by introducing the auxiliary primal problem:

(7.1) $$u_1(x^1) := u(x^1, 0_{d-1+N}), \quad x^1 \in \mathbb{R}_+,$$

and dualize the value function $u_1$ as follows. Our set of dual variables is defined as

(7.2) $$\mathcal{D}(y^1) = \bigg\{ (Y, \alpha) \in L^1(\Omega \times \mathbb{T}, \mathbb{R}^d_+) \times \mathbb{R}_+ : \forall x^1 \in \mathbb{R}_+,$$
$$\forall c \in \mathcal{C}_T((x^1, 0_{d-1+N}); K, R),$$
$$\mathbb{E}\bigg[\sum_{t \in \mathbb{T}} Y_t \cdot c_t - y^1 x^1\bigg] \leq \alpha \bigg\},$$
$$y^1 \in \mathbb{R}_+,$$

and we consider the dual problem

$$\tilde{u}_1(y^1) = \inf_{(Y,\alpha)\in\mathcal{D}(y^1)} \mathbb{E}\left[\sum_{t\in\mathbb{T}} \tilde{U}_t(Y_t) + \alpha\right], \qquad y^1 \in \mathbb{R}_+. \tag{7.3}$$

Recall that by convention $L^1(\Omega \times \mathbb{T}, \mathbb{R}^d_+) = L^1(\Omega \times \mathbb{T}, \mathbb{R}^d_+; \mathbb{P})$.

REMARK 7.1. By Remark 6.2, we can find some $x = (x^1, 0_{d-1+N}) \in \mathbb{R}^{d+N}$ such that the constant consumption process $c$ defined by $c_t^i = 1$ for all $t \in \mathbb{T}$ and $i \leq d$ belongs to $\mathcal{C}_T^U(x; K, R)$. It then follows from the definition of $\mathcal{D}(y^1)$ that, for each $\alpha \in \mathbb{R}_+$, the set $\{Y : (Y, \alpha) \in \mathcal{D}(y^1)\}$ is bounded in $L^1(\Omega \times \mathbb{T}, \mathbb{R}^d_+)$.

The abstract duality relation can be stated as follows.

LEMMA 7.1. *Under the assumptions of Theorem 4.1, we have the duality relations:*

$$\tilde{u}_1(y^1) = \sup_{x^1 \in \mathbb{R}_+} [u_1(x^1) - x^1 y^1], \qquad y^1 \in \mathbb{R}_+, \tag{7.4}$$

$$u_1(x^1) = \inf_{y^1 \in \mathbb{R}_+} [\tilde{u}_1(x^1) - x^1 y^1], \qquad x^1 \in \mathbb{R}_+. \tag{7.5}$$

PROOF. We only establish (7.4). The other relation (7.5) follows from (7.4) and general bidual properties of Legendre transform (see, e.g., [12]).

By definitions of $\tilde{U}_t$ and $\mathcal{D}(y^1)$, we have for all $x^1, y^1 \in \mathbb{R}_+$, $c \in \mathcal{C}_T((x^1, 0_{d-1+N}); K, R)$ and $(Y, \alpha) \in \mathcal{D}(y^1)$:

$$\mathbb{E}\left[\sum_{t\in\mathbb{T}} U_t(c_t)\right] \leq \mathbb{E}\left[\sum_{t\in\mathbb{T}} \tilde{U}_t(Y_t) + Y_t \cdot c_t\right]$$
$$\leq \mathbb{E}\left[\sum_{t\in\mathbb{T}} \tilde{U}_t(Y_t) + \alpha\right] + x^1 y^1, \tag{7.6}$$

and so

$$w(y^1) := \sup_{x^1 \in \mathbb{R}_+} [u_1(x^1) - x^1 y^1] \leq \tilde{u}_1(y^1) \qquad \forall y^1 \in \mathbb{R}_+. \tag{7.7}$$

We now fix some $y^1 \in \mathbb{R}_+$. In order to prove (7.4), we can assume w.l.o.g. that $w(y^1) < \infty$.

1. For $n > 0$, we define $\mathcal{C}_n$ as

$$\mathcal{C}_n = \{c = (c_t)_t \in L^0(\mathbb{R}^d_+; \mathbb{F}) : |c_t| \leq n, t \in \mathbb{T}\}.$$

The sets $\mathcal{C}_n$ are compact for the weak topology $\sigma(L^\infty(\Omega \times \mathbb{T}, \mathbb{R}_+^d), L^1(\Omega \times \mathbb{T}, \mathbb{R}_+^d))$. Moreover, it is clear from its definition that $\mathcal{D}(y^1)$ is a closed convex subset of $L^1(\Omega \times \mathbb{T}, \mathbb{R}_+^d)$. We may then apply the min–max theorem to get

$$\sup_{c \in \mathcal{C}_n} \inf_{(Y,\alpha) \in \mathcal{D}(y^1)} \mathbb{E}\left[\sum_{t \in \mathbb{T}}(U_t(c_t) - Y_t \cdot c_t) + \alpha\right]$$

$$= \inf_{(Y,\alpha) \in \mathcal{D}(y^1)} \sup_{c \in \mathcal{C}_n} \mathbb{E}\left[\sum_{t \in \mathbb{T}}(U_t(c_t) - Y_t \cdot c_t) + \alpha\right].$$

By setting

$$\tilde{U}_t^n(y) = \sup_{c \in \mathbb{R}_+^d, |c| \leq n} [U_t(c) - c \cdot y], \qquad y \in \mathbb{R}_+^d,$$

we then deduce that

(7.8)
$$\sup_{c \in \mathcal{C}_n} \inf_{(Y,\alpha) \in \mathcal{D}(y^1)} \mathbb{E}\left[\sum_{t \in \mathbb{T}}(U_t(c_t) - Y_t \cdot c_t) + \alpha\right]$$

$$= \inf_{(Y,\alpha) \in \mathcal{D}(y^1)} \mathbb{E}\left[\sum_{t \in \mathbb{T}} \tilde{U}_t^n(Y_t) + \alpha\right] := \tilde{u}_1^n(y^1).$$

For later use, observe that

(7.9) $\tilde{U}_t^n(y) \geq \tilde{U}_t^n(z)$ if $z - y \in \mathbb{R}_+^d$ and $\tilde{U}_t^n \leq \tilde{U}_t^k$ if $k \geq n$.

2. For any $Z = (Z^F, Z^I) \in \mathcal{Z}_T(K, \mathbb{P})$, we have from the $\mathbb{P}$-martingale property of $Z^F$ and Theorem 6.2: $\forall x^1 \in \mathbb{R}_+, \forall c \in \mathcal{C}_T((x^1, 0_{d-1+N}); K, R)$,

$$\mathbb{E}\left[Z_T^F \cdot \sum_{t \in \mathbb{T}} c_t - Z_0^1 x^1\right] = \mathbb{E}\left[\sum_{t \in \mathbb{T}} Z_t^F \cdot c_t - Z_0^1 x^1\right]$$

$$\leq a(0_N; Z, \mathbb{P}).$$

It follows that the pairs $(Y, \alpha)$ defined by

(7.10) $$Y = \frac{y^1}{Z_0^1} Z^F, \qquad \alpha = \frac{y^1}{Z_0^1} a(0_N; Z, \mathbb{P})$$

belong to $\mathcal{D}(y^1)$. Here, we use the convention $0/0 = 0$ and we observe from Remark 4.1 and the martingale property of $Z^F$ that $Z_0^1 = 0$ implies $Z^F = 0$.

Now, for $x^1 \in \mathbb{R}_+$, let $c = (c_t) \in L^0(\mathbb{R}_+^d; \mathbb{F})$ be such that

$$\mathbb{E}\left[\sum_{t \in \mathbb{T}} Y_t \cdot c_t - y^1 x^1\right] \leq \alpha \qquad \forall (Y, \alpha) \in \mathcal{D}(y^1).$$

By taking $(Y, \alpha)$ in the form (7.10), we deduce that

$$\mathbb{E}\left[Z_T^F \cdot \sum_{t \in \mathbb{T}} c_t - Z_0^1 x^1\right] \leq a(0_N; Z, P) \qquad \forall Z \in \mathcal{Z}_T(K, \mathbb{P}).$$

By Theorem 6.2, this means $(c_t) \in \mathcal{C}_T((x^1, 0_{d-1+N}); K, R)$. Therefore, we have the duality relation between the sets $\mathcal{C}(x^1)$ and $\mathcal{D}(y^1)$ in the sense that, for any $x^1 \in \mathbb{R}_+$, an element $c = (c_t)$ in $L^0(\mathbb{R}_+^d; \mathbb{F})$ belongs to $\mathcal{C}_T((x^1, 0_{d-1+N}); K, R)$ if and only if

$$\mathbb{E}\left[\sum_{t \in \mathbb{T}} Y_t \cdot c_t - y^1 x^1\right] \leq \alpha \qquad \forall (Y, \alpha) \in \mathcal{D}(y^1).$$

It follows that

$$\sup_{(Y,\alpha) \in \mathcal{D}(y^1)} \mathbb{E}\left[\sum_{t \in \mathbb{T}} Y_t \cdot c_t - \alpha\right]$$
$$= \inf\{y^1 x^1 : x^1 \geq 0 \text{ s.t. } c \in \mathcal{C}_T((x^1, 0_{d-1+N}); K, R)\},$$

and therefore

(7.11)
$$\lim_{n \to \infty} \sup_{c \in \mathcal{C}_n} \inf_{(Y,\alpha) \in \mathcal{D}(y^1)} \mathbb{E}\left[\sum_{t \in \mathbb{T}} (U_t(c_t) - Y_t \cdot c_t) + \alpha\right]$$
$$= \sup_{x_1 \in \mathbb{R}_+} \sup_{c \in \mathcal{C}_T((x^1, 0_{d-1+N}); K, R)} \mathbb{E}\left[\sum_{t \in \mathbb{T}} U_t(c_t) - x^1 y^1\right] = w(y^1).$$

3. Identifying relations (7.8) and (7.11), we get

(7.12) $$\lim_{n \to \infty} \tilde{u}_1^n(y^1) = w(y^1),$$

and so we have to show that

(7.13) $$\lim_{n \to \infty} \tilde{u}_1^n(y^1) = \tilde{u}_1(y^1).$$

Let $(Y^n, \alpha^n)$ be a sequence in $\mathcal{D}(y^1)$ such that

$$\lim_{n \to \infty} \mathbb{E}\left[\sum_{t \in \mathbb{T}} \tilde{U}_t^n(Y_t^n) + \alpha^n\right] = \lim_{n \to \infty} \tilde{u}_1^n(y^1) = w(y^1).$$

By Lemma A1.1 in [3] applied on $L^0(\Omega \times \mathbb{T})$, there exists a sequence $(\hat{Y}^n) \in \text{conv}(Y^n, \ldots, Y^{n+1}, \ldots)$ which converges a.e. to a process $\hat{Y}$, taking possibly infinite values. Moreover, by convexity of $\tilde{U}_t^n$ and (7.9), we have, by (7.12),

(7.14) $$\liminf_{n \to \infty} \mathbb{E}\left[\sum_{t \in \mathbb{T}} \tilde{U}_t^n(\hat{Y}_t^n) + \hat{\alpha}^n\right] \leq \lim_{n \to \infty} \tilde{u}_1^n(y^1) = w(y^1),$$

where $\hat{\alpha}^n$ is constructed from $(\alpha^k)_{k\geq n}$ with the same convex combinations as $\hat{Y}^n$.

Let us now consider the sequence of nonincreasing convex functions

$$\varphi_t := (-V_t)^{-1} \qquad \text{on } \mathbb{R}_+,$$

where $(-V_t)^{-1}$ denotes the generalized inverse of the nondecreasing function $-V_t$, and $V_t$ is the convex map introduced in Remark 4.2. We then define

$$\ell : y \in \mathbb{R}^d \mapsto \min\left\{x \cdot y : x \in \mathbb{R}^d \text{ with } \|x\| := \sum_{i=1}^d |x^i| e^i = 1\right\}.$$

With this notation, we have $y - \ell(y)e \in \mathbb{R}^d_+$ for all $y \in \mathbb{R}^d$. Since $\varphi_t$ is nonincreasing and $\tilde{U}^n_t \leq \tilde{U}_t$, it follows from (7.9) that

$$\begin{aligned}
\mathbb{E}[\varphi_t(\tilde{U}^n_t(\hat{Y}^n_t)^-)] &\leq \mathbb{E}[\varphi_t(\tilde{U}^n_t(\ell(\hat{Y}^n_t)e)^-)] \\
&\leq \varphi_t(0) + \mathbb{E}[\ell(\hat{Y}^n_t)] \\
&\leq \varphi_t(0) + \mathbb{E}[X(e) \cdot \hat{Y}^n_t] \qquad \text{with } X(e) = (1/e^1, \ldots, 1/e^d).
\end{aligned}$$

By part 2 of Remark 4.1, we can find some $x(e) > 0$ such that $(x(e), 0_{d-1+N}) - X(e) \in K_0$. Then, $(X(e))_{t\in\mathbb{T}} \in \mathcal{C}_T(Tx(e)\mathbf{1}_1; K, R)$ where $\mathbf{1}_1 = (1, 0, \ldots, 0) \in \mathbb{R}^{d+N}$. It then follows from the above inequality and the definition of $\mathcal{D}(y^1)$ that

$$(7.15) \qquad \mathbb{E}\left[\sum_{t\in\mathbb{T}} \varphi_t(\tilde{U}^n_t(\hat{Y}^n_t)^-)\right] \leq \sum_{t\in\mathbb{T}} \varphi_t(0) + Tx(e)y^1 + \hat{\alpha}^n.$$

Now, by l'Hôpital's rule and Remark 4.2, $\varphi_t(r)/r$ goes to infinity when $r$ goes to infinity, and so there exists some positive $\bar{r}_t > 0$ such that $\varphi_t(r) \geq 2r$ for all $r \geq \bar{r}_t$. Hence, for all $n$,

$$\tilde{U}^n_t(\hat{Y}^n_t)^- \leq \bar{r}_t + \tfrac{1}{2}\varphi_t(\tilde{U}^n_t(\hat{Y}^n_t)^-),$$

and by (7.15)

$$\mathbb{E}\left[\sum_{t\in\mathbb{T}} \tilde{U}^n_t(\hat{Y}^n_t)^-\right] \leq \bar{C}(y^1) + \tfrac{1}{2}\hat{\alpha}^n,$$

where $\bar{C}(y^1) = \sum_{t\in\mathbb{T}} \bar{r}_t + \tfrac{1}{2}(\varphi_t(0) + x(e)y^1)$. We then deduce that

$$\tfrac{1}{2}\hat{\alpha}^n - \bar{C}(y^1) \leq \mathbb{E}\left[\sum_{t\in\mathbb{T}} \tilde{U}^n_t(\hat{Y}^n_t)\right] + \hat{\alpha}^n$$

so that by (7.14), after possibly passing to a subsequence, for $n$ large enough

$$\tfrac{1}{2}\hat{\alpha}^n \leq w(y^1) + 1 + \bar{C}(y^1) < \infty,$$

which proves that the sequence $(\hat{\alpha}_n)$ is bounded. After possibly passing to a subsequence, we can then assume that it converges to some $\hat{\alpha} \in \mathbb{R}_+$. It then follows from (7.15) and the de la Vallée Poussin theorem that the sequence $(\tilde{U}_t^n(\hat{Y}^n)^-)$ is uniformly integrable. Since $\tilde{U}_t^n$ converges to $\tilde{U}_t$ uniformly on compact sets, it then follows from Fatou's lemma that

$$\liminf_{n \to \infty} \mathbb{E}\left[\sum_{t \in \mathbb{T}} \tilde{U}_t^n(\hat{Y}_t^n) + \hat{\alpha}^n\right] \geq \mathbb{E}\left[\sum_{t \in \mathbb{T}} \tilde{U}_t(\hat{Y}_t) + \hat{\alpha}\right] \geq \tilde{u}_1(y^1).$$

Since we obviously have $\tilde{u}_1^n(y^1) \leq \tilde{u}_1(y^1)$, inequality (7.14) implies (7.13), that is, (7.4), see (7.7). □

REMARK 7.2. 1. Assume that for some $\hat{x}^1 > 0$, $u_1(\hat{x}^1) < \infty$. Then, by the duality relation (7.5) in Lemma 7.1, there exists some $y^1 \in \mathbb{R}_+$ such that $\tilde{u}_1(y^1) < \infty$. Hence, for this $y^1$, there exists some $(Y, \alpha) \in \mathcal{D}(y^1)$ such that $\sum_{t \in \mathbb{T}} \tilde{U}_t(Y_t) \in L^1(\mathbb{P})$. In view of (7.6), this implies that $u_1(\tilde{x}^1) < \infty$ for all $\tilde{x}^1 \geq 0$.

2. Fix $x \in \text{int}(K_0)$. Then, by parts 2 and 3 of Remark 4.1, there exists some $\hat{x} = (\hat{x}^1, 0_{d-1+N})$ with $\hat{x}^1 > 0$ such that $x - (\hat{x}^1, 0_{d-1+N}) \in K_0$ and $\mathcal{C}_T^U(\hat{x}; K, R) \neq \varnothing$. Since $\mathcal{C}_T(\hat{x}; K, R) \subset \mathcal{C}_T(x; K, R)$, the finiteness of $u(x)$ implies the finiteness of $u_1(\hat{x}^1) = u(\hat{x})$.

3. Finally, let $x \in \mathbb{R}^d \times \mathbb{R}_+^N$ be such that $\mathcal{C}_T^U(x; K, R) \neq \varnothing$. Then, by Remark 6.2, there is some $\tilde{x}^1 > 0$ such that $u(\tilde{x}^1, 0_{d-1+N}) \geq u(x)$.

4. Combining part 2 with part 1 and then part 3 proves (i) of Theorem 4.1.

We go on preparing the proof of Theorem 4.1 with two more lemmas.

LEMMA 7.2. *Let the conditions* (4.1), (4.3) *and* (4.5) *hold. Let* $y^1 \in \mathbb{R}_+$ *and* $(Y, \alpha) \in \mathcal{D}(y^1)$ *be such that*

$$\left(\sum_{t \in \mathbb{T}} \tilde{U}_t(Y_t)\right) \in L^1(\mathbb{P}).$$

*Then,*

$$\left(\sum_{t \in \mathbb{T}} \tilde{U}_t(\lambda Y_t)\right) \in L^1(\mathbb{P}) \qquad \textit{for all } \lambda \in (0, 1].$$

PROOF. 1. First observe that $\tilde{U}_t(\lambda Y_t)^- \in L^1(\mathbb{P})$ for each $t \in \mathbb{T}$ and $\lambda \in (0, 1]$. Indeed, given $x^F \in \text{int}(\mathbb{R}_+^d)$, we have by definition of $\tilde{U}_t$

$$U_t(x^F) \leq \tilde{U}_t(\lambda Y_t) + \lambda Y_t \cdot x^F.$$

Since $Y_t \in L^1(\mathbb{P})$ and $U_t(x^F)$ is finite, this implies that $\tilde{U}_t(\lambda Y_t)^- \in L^1(\mathbb{P})$.

2. From part 1, it suffices to show that $\tilde{U}_t(Y_t)^+ \in L^1(\mathbb{P})$ implies that $\tilde{U}_t(\lambda Y_t)^+ \in L^1(\mathbb{P})$ for all $\lambda \in (0,1]$. Fix $\lambda \in (0,1]$; by assumption (4.5) we have

$$\tilde{U}_t(\lambda Y_t)^+ \leq C_t^\lambda(1 + \tilde{U}_t(Y_t)^+),$$

which, by part 1, shows that $\tilde{U}_t(\lambda Y_t)^+ \in L^1(\mathbb{P})$ and concludes the proof. $\square$

LEMMA 7.3. *Fix $(K,R) \in \mathcal{K} \times \mathcal{R}$ such that* (R1)–(R3) *and $NA^r(K,R)$ hold. Fix $x \in \mathbb{R}^d \times \mathbb{R}_+^N$, and let $(c_n)_{n\geq 1}$ be a sequence in $\mathcal{C}_T(x;K,R)$. Then, there is a sequence $(\tilde{c}_n)_{n\geq 1}$ such that $\tilde{c}_n \in \mathrm{conv}(c_k, k \geq n)$, for each $n \geq 1$, which converges $\mathbb{P}$-a.s. to some $\tilde{c} \in \mathcal{C}_T(x;K,R)$.*

PROOF. We set $x = (x^F, x^I)$. Since $c_n^i \geq 0$ for each $1 \leq i \leq d+N$, we deduce from Lemma A1.1 in [3] that there is a sequence $(\tilde{c}_n)_{n\geq 1}$ such that $\tilde{c}_n \in \mathrm{conv}(c_k, k \geq n)$, for each $n \geq 1$, which converges $\mathbb{P}$-a.s. to some $\tilde{c} \in L^0([0,\infty];\mathbb{F})$. By Lemma 2.1, $\tilde{c}_n \in \mathcal{C}_T(x;K,R)$ for each $n \geq 1$. Since, by Remark 5.2, $\mathcal{C}_T(x;K,R)$ is closed, it suffices to show that $\|\sum_{t\in\mathbb{T}} \tilde{c}_t\| < \infty$. To see this, recall from Theorem 6.2 and Lemma 6.1 that there is some $Z = (Z^F, Z^I) \in \mathcal{Z}_T(K,\mathbb{P})$ such that

$$\mathbb{E}\left[Z_T^F \cdot \left(\sum_{t\in\mathbb{T}}(c_n)_t\right)\right] \leq Z_0^F \cdot x^F + a(x^I; Z, \mathbb{P}) < \infty, \qquad n \geq 1.$$

By Remark 6.1, we have $Z_T^i > 0$ $\mathbb{P}$-a.s. for all $1 \leq i \leq d$; sending $n$ to $\infty$ and using Fatou's lemma then leads to the required result. $\square$

We can now conclude the proof of Theorem 4.1.

PROOF OF THEOREM 4.1. Item (i) has already been proved in Remark 7.2. We prove (ii).

1. Let $(c^n)_{n\geq 1}$ be a sequence in $\mathcal{C}_T(x;K,R)$ such that

$$u(x) = \lim_{n\to\infty} \mathbb{E}\left[\sum_{t\in\mathbb{T}} U_t(c_t^n)\right].$$

Since $U_t$ is convex, it follows from Lemma 7.3 that, after possibly passing to convex combinations, we can assume that $c^n$ converges $\mathbb{P}$-a.s. to some $c^* \in \mathcal{C}_T(x;K,R)$. We shall prove in part 2 that

$$(7.16) \qquad \left\{\left(\sum_{t\in\mathbb{T}} U_t(c_t^n)\right)^+\right\}_{n\geq 1} \text{ is uniformly integrable.}$$

Then, using Fatou's lemma and the continuity of $U_t$, we obtain

$$u(x) = \lim_{n\to\infty} \mathbb{E}\left[\sum_{t\in\mathbb{T}} U_t(c_t^n)\right] \leq \mathbb{E}\left[\limsup_{n\to\infty} \sum_{t\in\mathbb{T}} U_t(c_t^n)\right] = \mathbb{E}\left[\sum_{t\in\mathbb{T}} U_t(c_t^*)\right].$$

2. To prove (7.16), we assume to the contrary that the sequence is not uniformly integrable and work toward a contradiction. If (7.16) does not hold, then, after possibly passing to a subsequence, we can find some $\delta > 0$ and a sequence $(A_k^n)_{k,n}$ such that, for each $n \geq 1$, $(A_k^n)_{k=1}^n$ forms a disjoint partition of $\Omega$ such that

$$(7.17) \qquad \mathbb{E}\left[\left(\sum_{t \in \mathbb{T}} U_t(c_t^k)\right)^+ \mathbb{1}_{A_k^n}\right] \geq \delta, \qquad 1 \leq k \leq n, n \geq 1.$$

By possibly adding a constant to the $U_t$'s, we can assume that there is some $r \in \mathbb{R}_+^d$ such that $\min_{t \in \mathbb{T}} U_t(r) \geq 0$.

Now, by Remark 7.2, there exists some $y^1 \in \mathbb{R}_+$ such that $\tilde{u}_1(y^1) < \infty$. Hence, for this $y^1$, there exists some $(Y, \alpha) \in \mathcal{D}(y^1)$ such that

$$\left(\sum_{t \in \mathbb{T}} \tilde{U}_t(Y_t)\right) \in L^1(\mathbb{P}).$$

Then, by Lemma 7.2,

$$\left(\sum_{t \in \mathbb{T}} \tilde{U}_t(\lambda Y_t)\right) \in L^1(\mathbb{P}) \qquad \text{for all } \lambda \in (0, 1].$$

Observe from Remark 6.2 that we can find some $\hat{x}^1 > 0$ such that the process $r + c^k$ belongs to $\mathcal{C}_T((\hat{x}^1, 0_{d-1+N}); R, K)$ for all $k \geq 1$. It then follows, by definitions of $\tilde{U}_t$ and $\mathcal{D}(y^1)$, that for each $\lambda > 0$ and $n \geq 1$

$$\sum_{k=1}^n \mathbb{E}\left[\sum_{t \in \mathbb{T}} U_t(r + c_t^k) \mathbb{1}_{A_k^n}\right] \leq \mathbb{E}\left[\sum_{t \in \mathbb{T}} \tilde{U}_t(\lambda Y_t)\right] + \lambda \sum_{k=1}^n \mathbb{E}\left[Y_T \cdot \left(\sum_{t \in \mathbb{T}} r + c_t^k\right) \mathbb{1}_{A_k^n}\right]$$

$$\leq \mathbb{E}\left[\sum_{t \in \mathbb{T}} \tilde{U}_t(\lambda Y_t)\right] + n\lambda(y^1 \hat{x}^1 + \alpha).$$

Since $U_t$ is $\mathbb{R}^d$-nondecreasing [see (4.2)], we have $U_t(r + c_t^k) \geq U_t(c_t^k)^+$. It then follows from (7.17) that

$$n\delta \leq \mathbb{E}\left[\sum_{t \in \mathbb{T}} \tilde{U}_t(\lambda Y_t)\right] + n\lambda(y^1 \hat{x}^1 + \alpha) \qquad \text{for all } n \geq 1 \text{ and } \lambda > 0.$$

Dividing by $n \geq 1$ and sending $n$ to $\infty$ in the above inequality, we obtain

$$(7.18) \qquad\qquad \delta \leq \lambda(y^1 \hat{x}^1 + \alpha) \qquad \text{for all } \lambda > 0.$$

Sending $\lambda$ to 0 then leads to the required contradiction since $\delta > 0$. $\square$

REMARK 7.3. 1. Since $\mathcal{D}(\lambda y^1) = \lambda \mathcal{D}(y^1)$ for all $\lambda \geq 1$, the above proof goes through if we replace the assumption (4.5) by

$$\tilde{u}^1(y^1) < \infty \qquad \text{for all } y^1 > 0. \tag{7.19}$$

Moreover, as explained above, it follows from Remark 6.2 that $u(x) < \infty$ whenever $\tilde{u}^1(y^1) < \infty$ for some $y^1 \geq 0$. Hence, if (7.19) holds, then the assumption $u(x) < \infty$ for some $x \in \text{int}(K_0)$ can be dropped too.

2. Since $\tilde{u}^1$ is nonincreasing, it follows from Lemma 7.2 that (7.19) is implied by (4.5) and the condition $u(x) < \infty$ for some $x \in \text{int}(K_0)$.

Laboratoire de Probabilités
  et Modèles Aléatoires
CNRS, UMR 7599
Université Paris 6
and
CREST
France
e-mail: bouchard@ccr.jussieu.fr

Laboratoire de Probabilités
  et Modèles Aléatoires
CNRS, UMR 7599
Université Paris 7
and
CREST
France
e-mail: pham@math.jussieu.fr